\documentclass{amsart}
\usepackage{amsthm,amssymb,latexsym,amsxtra}

\usepackage[cmtip,all]{xy}

\newcommand{\lto}{\longrightarrow}
 
\DeclareMathOperator{\rank}{rank}
\DeclareMathOperator{\grade}{grade}
\DeclareMathOperator{\depth}{depth}
\DeclareMathOperator{\height}{height}
\DeclareMathOperator{\codim}{codim}
\DeclareMathOperator{\Spec}{Spec}
\DeclareMathOperator{\projdim}{proj\; dim}

\theoremstyle{plain}

\newtheorem{thm}{Theorem}[section]
\newtheorem{cor}[thm]{Corollary}
\newtheorem{lem}[thm]{Lemma}
\newtheorem{pro}[thm]{Proposition}

\theoremstyle{definition}
\newtheorem{dfn}[thm]{Definition}

\newtheorem{block}[thm]{}
\newtheorem*{dfn0}{Definition}
\newtheorem*{rem0}{Remark}

\theoremstyle{remark}
\newtheorem*{notation}{Notation}
\newtheorem*{acknowledgments}{Acknowledgments}

\begin{document}
\title[Unirationality of Hurwitz spaces]{Unirationality of Hurwitz spaces\\ 
of coverings of degree $\leq 5$}
\author{V. Kanev}
\date{}
\address{V. Kanev\\Dipartamento di Matematica\\
         Universit\`{a} di Palermo\\Via Archirafi, 34\\
         90123 Palermo\\Italy}
                                \email{vassil.kanev@unipa.it}
                                

\dedicatory{To the memory of my friend and teacher,\\ Andrey N. Todorov}

\begin{abstract}
Let $Y$ be a smooth, projective curve of genus $g\geq 1$ over the complex numbers. 
Let $\mathcal{H}^0_{d,A}(Y)$
be the Hurwitz space which parameterizes  
equivalence classes of
coverings $\pi :X \to Y$ of
degree $d$ simply branched in $n=2e$ points, such that 
the monodromy group is $S_d$ and 
$\det(\pi_{*}\mathcal{O}_X/\mathcal{O}_Y)^{\vee}$ is 
isomorphic to a fixed line bundle $A$ of degree $e$. We prove that, 
when $d=3, 4$ or $5$ and $n$ is sufficiently large (precise bounds are given),
the Hurwitz space $\mathcal{H}^0_{d,A}(Y)$ is unirational. If in addition $(e,2)=1$ (when $d=3$),
$(e,6)=1$ (when $d=4$) and $(e,10)=1$ (when $d=5$), then  $\mathcal{H}^0_{d,A}(Y)$ is
rational.
\end{abstract}

\maketitle

\def\thefootnote{}
\footnote{The author participated in the research project
PRIN 2008-12 `\emph{Geometria delle variet\`{a} algebriche e dei loro spazi di moduli}'}

\section*{Introduction}
The Hurwitz space $\mathcal{H}_{d,n}(Y)$ parameterizes classes $[X\to Y]$ of 
equivalent coverings of degree $d$ of a smooth, projective curve $Y$, 
simply branched in $n$ points, with smooth, irreducible $X$. The irreducibility of 
$\mathcal{H}_{d,n}(\mathbb{P}^{1})$ was proved by Hurwitz in 
\cite{Hur}, based on earlier work by Clebsch and L\"{u}roth. The 
unirationality of $\mathcal{H}_{d,n}(\mathbb{P}^{1})$, when $d\leq 3$ 
is classically known (cf. \cite{Mi}). Arbarello and Cornalba proved in \cite{AC}
the unirationality of $\mathcal{H}_{d,n}(\mathbb{P}^{1})$ when $d=4$ 
or 5.\, F.-O. Schreyer gave in \cite{Sch} another proof of this result. 
\par
In studying Hurwitz spaces when $g(Y)\geq 1$, it is natural to restrict 
ourselves to 
simple coverings of degree $d$ whose monodromy group is $S_d$. We denote the
corresponding Hurwitz space by $\mathcal{H}^0_{d,n}(Y)$. This is the principal case, 
since those with smaller monodromy group are reduced to coverings of smaller 
degree via an \'{e}tale covering of $Y$. 
There is  a 
canonical morphism $\mathcal{H}^0_{d,n}(Y)\to Pic^{\frac{n}{2}}Y$.
If $[\pi:X\to Y]\in \mathcal{H}^0_{d,n}(Y)$ and $\pi_*\mathcal{O}_{X}/\mathcal{O}_{Y}$ is its
\emph{Tschirnhausen module} (as named in \cite{Mi}), then the morphism sends
$[\pi:X\to Y]$ to $\det(\pi_*\mathcal{O}_{X}/\mathcal{O}_{Y})^{-1}$.
 The problem of 
unirationality may be posed for the fibers of this morphism. Given 
$A\in Pic^{\frac{n}{2}}Y$ we denote the fiber over $A$ by 
$\mathcal{H}^0_{d,A}(Y)$. 
\par
We study the problem of irreducibility and 
unirationality/rationality of $\mathcal{H}_{d,A}^0(Y)$. The case $d=2$ is trivial, 
since $\mathcal{H}_{2,A}^0(Y)=\mathcal{H}_{2,A}(Y)$ and this variety 
is isomorphic to $|A^{2}|\setminus \Delta$, where $\Delta$ is the 
codimension one subvariety of $|A^{2}|$, parameterizing divisors with 
multiplicities (see e.g. \cite{Wa}). The author studied in \cite{K1} 
and \cite{K2} the Hurwitz spaces $\mathcal{H}_{3,A}^0(Y)$ and 
$\mathcal{H}_{4,A}^0(Y)$ when $g(Y)=1$ and $\deg A > 0$ and proved their 
irreducibility and unirationality. These results were applied, in the 
same papers, to proving the unirationality of the Siegel modular 
varieties $\mathcal{A}_{3}(1,1,d)$ and $\mathcal{A}_{3}(1,d,d)$ when 
$d\leq 4$. 
The case $\deg A=3$ is the relevant one for these applications.
\par
We prove the following theorems.
\begin{thm}\label{s0.3}
Let $Y$ be a smooth, irreducible, projective curve of genus $g\geq 1$. 
Let $d=3,\, 4$ or $5$. Suppose the integer $e$ satisfies the 
inequality
\[
e >
\begin{cases}
4(g-1)+4          & \text{if $d=3$}\\
12(g-1)+6         & \text{if $d=4$}\\
40(g-1)+20        & \text{if $d=5$}.
\end{cases}
\]
If $A\in Pic^{e}Y$,
let $\mathcal{H}_{d,A}^0(Y)$ be the Hurwitz space,
parameterizing equivalence classes of simply ramified coverings of $Y$, 
branched in $n=2e$ points, with full monodromy group $S_d$, and 
determinants of the Tschirnhausen modules isomorphic to $A^{-1}$. Then 
$\mathcal{H}_{d,A}^0(Y)$ is irreducible and unirational. 
\end{thm}
\begin{thm}\label{s0.100}
Under the hypothesis of Theorem~\ref{s0.3} assume moreover that 
\begin{alignat*}{2}
(e,2)&=1 &\quad&\text{if\; $d=3$} \\
(e,6)&=1 &&\text{if\; $d=4$} \\
(e,10)&=1 &&\text{if\; $d=5$}. 
\end{alignat*}
Let $A\in Pic^{e}Y$. Then $\mathcal{H}_{d,A}^0(Y)$ is a rational 
variety.
\end{thm}
The proof of these results is based on the classification of the 
Gorenstein coverings of degree $\leq 5$, by means of vector bundles 
over the base, due to Miranda, Casnati and Ekedahl \cite{Mi}, 
\cite{CE}, \cite{C2}. The unirationality/rationality of the Hurwitz
spaces is eventually a consequence of  known results on unirationality/rationality of 
moduli spaces of vector bundles  over curves with fixed determinants.
Here is an outline of the content of the paper by sections. 
\par
Section~\ref{s1} contains
some preliminaries on vector bundles. 
\par
Section~\ref{s5} is devoted to the explicit
description of the Gorenstein coverings of degree $d=3, 4$ or $5$. We recall and complement some results from the papers \cite{Mi}, \cite{CE} and \cite{C2}. In the proofs of Theorem~\ref{s0.3}
and Theorem~\ref{s0.100} we use certain statements about  uniqueness of the representation of a Gorenstein covering  by means of the vector bundle data
on the base. When $d=5$
this is related with a problem discussed in  \cite{BuE} p.457 Remark~2. Namely, given the
Buchsbaum-Eisenbud resolution of a Gorenstein ideal of grade 3, to what extent is
the skew-symmetric matrix of odd degree uniquely
determined  by the Pfaffian ideal? We focus mainly on the more difficult case $d=5$
in this section
and prove in Lemma~\ref{s5.84}, Lemma~\ref{s5.87a}
and Proposition~\ref{s5.89b} the uniqueness statements we need. For reader's convenience we give in Proposition~\ref{BE3} a  proof of the Buchsbaum-Eisenbud resolution as we use it.
\par
In Section~\ref{s2} we construct families of coverings with rational parameter varieties, which eventually will dominate the Hurwitz spaces. Here the main result is Lemma~\ref{s2.35}.
\par
In Section~\ref{s3} the proofs of Theorem~\ref{s0.3} and Theorem~\ref{s0.100} are given.
\par
In Appendix~A we give a proof of a result due to Dolgachev and Libgober that we use. The reason of including this appendix is a comment in \cite{Sh} p.337 regarding its validity.  We give a detailed proof along the sketch in \cite{DL} p.9.

\begin{notation}
We assume the base field $k=\mathbb{C}$, unless otherwise specified. 
A scheme is always supposed to be separated of finite type.
  Points of a
scheme are always closed points. A variety is a reduced scheme.
A morphism of schemes $\pi : X\to Y$ is called covering if it is finite, 
flat and surjective. Two coverings $\pi_1 : X_1\to Y$ and $\pi_2 : 
X_2\to Y$ 
are equivalent
if there exists an isomorphism $f:X_1\to X_2$ such that $\pi_1=\pi_2\circ f$.
A locally free sheaf is always supposed to be 
coherent. 
 We make 
distinction between locally free sheaves and 
vector bundles and we denote differently their projectivizations. If $E$ is a 
locally free sheaf on $Y$ and if $\mathbb{E}$ is the corresponding vector 
bundle, i.e. $E\cong \mathcal{O}_Y(\mathbb{E})$, then $\mathbf{P}(E):=
\mathbf{Proj}(S(E)) \cong \mathbb{P}(\mathbb{E}^{\vee})$. We denote by 
$ad'(E)$ the sheaf of endomorphisms of $E$ with trace 0.
\end{notation}
\section{Preliminaries}\label{s1}
In this section we collect some facts about vector bundles we need. Contrary
to the rest of the paper we shall make here  the customary identification 
between
vector bundles and the associated locally free sheaves of sections. 
We include some 
proofs for which  we could not find references. 
\begin{block}\label{s1.13}
Let $E$ and $F$ be fixed vector bundles on a smooth, projective curve 
$Y$. It is well-known that the extensions $0\to E\to W\to F\to 0$ are 
classified by $H^1(Y,Hom(F,E))$ (see \cite{At2}). There is a canonical 
way to give a structure of an algebraic family on the set of 
extensions (see \cite{NS} Lemma~3.1), which we now recall. Let 
$H=H^1(Y,Hom(F,E))$. Let $e_{1},\ldots,e_{k}$ be a basis of $H$ and 
let $x_{i}, i=1,\ldots,k$ be the dual basis of $H^{*}$. One considers 
the element 
\begin{equation}\label{es1.13a}
\omega = \sum_{i=1}^{k}x_{i}\otimes e_{i} \in H^{*}\otimes 
H^1(Y,Hom(F,E))
\end{equation}
The projection map $p_Y:H\times Y\to Y$ is an affine morphism,
so the right-hand space of \eqref{es1.13a} is canonically embedded in 
$H^1(H\times Y, Hom(p_{Y}^{*}F,p_{Y}^{*}E))$. The extension of vector 
bundles over $H\times Y$ corresponding to $\omega$
\begin{equation}\label{es1.13}
0\to p_{Y}^{*}E\to V\to p_{Y}^{*}F\to 0
\end{equation}
has the property that for every $h\in H$ the restriction of 
\eqref{es1.13} on $\{h\}\times Y$ is an extension of $F$ with kernel 
$E$ whose class equals $h$.
\end{block}
The following lemma may be deduced from \cite{N0} Theorem~A on p.263.
\begin{lem}\label{s1.14}
Let $a\in H^1(Y,Hom(F,E))$ be the class of 
\begin{equation}\label{es1.14}
0\to E\overset{i}{\to} 
W\overset{p}{\to} F\to 0.
\end{equation}
Let $\kappa :T_{a}H\to H^1(Y,Hom(W,W))$ be 
the Kodaira-Spencer map associated with the family of vector bundles 
$V\to H\times Y$ (see \cite{KS}). Then $\kappa$ equals the canonical map $p^{*}\otimes 
i_{*}:H^1(Y,Hom(F,E))\to H^1(Y,Hom(W,W))$. 
\end{lem}
\begin{dfn}\label{s1.18}
Let $Y$ be an elliptic curve. Let $V$ be a vector bundle over $Y$ of 
rank $r$ and degree $e$. Let $h=g.c.d.(r,e)$ and let $r=r'h, e=e'h$. 
We say $V$ is \emph{regular polystable} if 
\[
V \cong V_{1}\oplus \cdots \oplus V_{h}
\]
where every $V_{i}$ is indecomposable of rank $r'$ and degree $e'$, 
and $V_{i}\ncong V_{j}$ if $i\neq j$. 
\end{dfn}
A vector bundle over an elliptic curve is stable if and only if it is 
indecomposable and its rank and degree are coprime (see e.g. 
\cite{Br} p.87). So, if $g.c.d.(r,e)>1$, no stable bundle of rank 
$r$ and degree $e$ over $Y$ exists. The notion of regular polystability 
is a replacement for the stability of vector bundles over elliptic 
curves. Every regular polystable vector bundle is semistable, since it 
is a direct sum of stable vector bundles of the same slope. 
\renewcommand{\theenumi}{\roman{enumi}}
\begin{pro}\label{s1.18aa}
Let $Y$ be an elliptic curve and let $V$ be a regular polystable 
bundle with decomposition $V=V_{1}\oplus \cdots \oplus V_{h}$ as in 
Definition~\ref{s1.18}. Then there exists a complex manifold $M$ and a 
holomorphic family $\mathcal{V}\to M\times Y$ of regular polystable 
bundles satisfying the following properties:
\begin{enumerate}
\item
there exists a point $m_{0}\in M$ such that $\mathcal{V}_{m_{0}}\cong 
V$;
\item
if $m_{1},m_{2}\in M$ and $m_{1}\neq m_{2}$, then 
$\mathcal{V}_{m_{1}}$ and $\mathcal{V}_{m_{2}}$ are not isomorphic;
\item
the Kodaira-Spencer map $\kappa:T_{m_{0}}M\to H^1(Y,End\, V)$ is an 
isomorphism;
\item
the family $\mathcal{V}\to M\times Y$ is complete at $m_{0}$. 
\end{enumerate}
\end{pro}
\begin{proof}
Let $J=Pic^{0}Y$ and let $\mathcal{L}\to J\times Y$ be the Poincar\'{e} 
line bundle of degree zero line bundles on $Y$. Consider the vector 
bundle $\mathcal{W}\to J^{h}\times Y$ defined as 
\[
\mathcal{W} = \oplus_{i=1}^{h}p_{Y}^{*}V_{i}\otimes (p_{i}\times 
id)^{*}\mathcal{L}
\]
where $p_{i}:J^{h}\to J$ is the projection onto the $i$-th factor. If $E$ is an 
indecomposable bundle of rank $n$ and degree $d$, such that $(n,d)=1$, 
then $E\otimes L\cong E$ if and only if $L$ is an $n$-torsion point, 
$L\in J_{n}$ (see \cite{At} p.434). Choose a disk $\Delta$ in $J$ 
which contains $0$ and satisfies $\Delta\cap J_{r'}=\{0\}$. Let 
$M=\Delta \times \cdots \times \Delta$ ($h$ times) and let 
$\mathcal{V}=\mathcal{W}|_{M\times Y}$. Then the family 
$\mathcal{V}\to M\times Y$ satisfies properties (i) and (ii) with 
$m_{0}=(0,\ldots,0)$. Statement (iv) follows from Statement (iii) 
according to \cite{NS} Lemma~2.1. So, it remains to verify (iii). If 
$i\neq j$, then $h^1(V_{i}^{*}\otimes V_{j}) = h^0(V_{i}^{*}\otimes 
V_{j}) = 0$ since $V_{i}\ncong V_{j}$. Hence $H^1(Y,End\, V) \cong 
\oplus_{i=1}^{n}H^1(Y,End\, V_{i})$. 

Let $E$ be an indecomposable 
bundle of rank $n$ and degree $d$ such that $(n,d)=1$. Atiyah proved 
that $End\, E\cong \oplus_{i=1}^{n^{2}}L_{i}$, where $\{L_{i}\}_{i}$ is 
the set of $n$-torsion points of $Pic^{0}Y$. Hence the trace 
homomorphism 
\[
Tr^* : H^1(Y,End\, E) \lto H^1(Y,\mathcal{O}_{Y})
\]
is an isomorphism. Consider the family $p_{Y}^{*}E\otimes 
\mathcal{L}\to J\times Y$ of vector bundles over $Y$. Its determinant \linebreak
$p_{Y}^{*}(\det E)\otimes \mathcal{L}^{n}$ is a family of line bundles 
of degree $d$ parameterized by $J$. The two families define 
Kodaira-Spencer maps, which fit into the following commutative diagram:
\[
\xymatrix{
&T_0J\ar[dl]_{\lambda}\ar[dr]^{\mu}\\
H^1(Y,End\, E)\ar[rr]^-{Tr^*}&&H^1(Y,\mathcal{O}_C)
}
\]
Since $\mu$ is a multiplication by $n$ one obtains that $\lambda:T_0J\to 
H^1(Y,End\, E)$ 
is an isomorphism. Applying this to every $V_i$ and taking the direct sum one 
concludes that
$
\kappa : T_0J^n \to H^1(Y,End\, E)
$
is an isomorphism.
\end{proof}
We need a sharpening of Proposition~2.6 of \cite{NR1}. We denote by $I^m$ the trivial vector bundle of rank $m$.
\renewcommand{\theenumi}{\roman{enumi}}
\begin{pro}\label{s1.18a}
Let $Y$ be a smooth, projective curve of genus $g\geq 1$. Let $W\to 
T\times Y$ be a family of vector bundles of rank $r$ and degree $e$. 
Then there exists an irreducible, nonsingular variety $S$ and a family 
of vector bundles $V\to S\times Y$ with the following properties.
\begin{enumerate}
\item
Every $W_{t}$ is isomorphic to $V_{s}$ for some $s\in S$.
\item
Every semistable bundle of rank $r$ and degree $e$ is isomorphic to 
$V_{s}$ for some $s\in S$.
\item
The variety $S$ is a vector bundle over the Jacobian $J$ of $Y$.
\item
For every $A\in Pic^{e}Y$ the variety $S_A = \{s\in S|\det V_{s}\cong A\}$ 
is isomorphic to a fiber of the vector bundle $S\to J$.
\item
If $g\geq 2$,  the open set $\{s\in S|V_{s}\; \text{is stable}\}$ is dense in 
$S$, and for every $A\in Pic^eY$,  the 
open set $\{s\in S_{A}| V_{s}\; \text{is stable}\}$ 
is dense in $S_{A}$. If $g=1$ the same statements hold replacing ``stable''
 with ``regular polystable''.
\item
For every $s\in S$, such that either $V_{s}$ is semistable, or $V_{s}$ 
is isomorphic to some $W_{t}$, the Kodaira-Spencer map $\kappa:T_{s}S 
\to H^1(Y,End\, V_{s})$ is epimorphic.
\item
For every $A\in Pic^{e}Y$ and every $s\in S_{A}$, such that either 
$V_{s}$ is semistable, or $V_{s}$ is isomorphic to some $W_{t}$, the 
Kodaira-Spencer map $\kappa:T_{s}S_{A}\to H^1(Y,End\, V_{s})$ has 
image $H^1(Y,ad'(V_{s}))$. 
\end{enumerate}
\end{pro}
\begin{proof}
If one proves the statements of the proposition for some pair $(r,e)$, 
then tensoring all bundles by a fixed line bundle $L$, one proves the 
statements for the pair $(r,e+r\deg L)$. Up to tensoring by a fixed line 
bundle of sufficiently large degree we may assume that the bundles of 
the family $\{W _{t}\}$ are globally generated and $H^1(Y,W_{t})=0$. 
Furthermore we may assume that $\mu = \frac{e}{r}> 2g-1$, so every 
semistable bundle $G$ of rank $r$ and degree $e$ is globally generated 
and $H^1(Y,G)=0$. By \cite{At} Theorem~2 every $W_{t}$, as well as 
every semistable $G$ with $rk(G)=r,\deg G=e$ contains a trivial bundle 
of rank $r-1$. Statements (i) through (v) are proved in \cite{NR1} 
Proposition~2.6 in the case $g\geq 2$. The same proof applies also to 
the case $g=1$, using the fact that regular polystability is a Zariski 
open condition in families of vector bundles over elliptic curves (see 
e.g. \cite{K1} Appendix~B). We notice that by construction the 
variety $S_A=\delta^{-1}(A)$ is isomorphic to $H^1(Y,Hom(A,I^{r-1}))$.
Let us first prove (vii). We have an exact sequence
\begin{equation}\label{es1.19a}
0\to I^{r-1}\overset{i}{\to} V_{s}\overset{p}{\to} A\to 0.
\end{equation}
According to Lemma~\ref{s1.14} the Kodaira-Spencer map 
$\kappa:T_{s}S_{A}\to H^1(Y,Hom(V_{s},V_{s}))$ equals
\[
p^{*}\otimes i_{*}:H^1(Y,Hom(A,I^{r-1}))\to H^1(Y,Hom(V_{s},V_{s})).
\]
From the exact sequence
\[
0\to Hom(A,V_{s})\overset{p^{*}}{\to} Hom(V_{s},V_{s})\to 
Hom(I^{r-1},V_{s})\to 0
\]
we obtain
\[
H^1(Y,Hom(A,V_{s}))\to H^1(Y,Hom(V_{s},V_{s}))\to 
H^1(Y,Hom(I^{r-1},V_{s}))\to 0.
\]
We have $H^1(Y,V_{s})=0$, so $H^1(Hom(A,V_{s}))\to 
H^1(Hom(V_{s},V_{s}))$ is surjective. The exact sequence 
\eqref{es1.19a} yields a commutative diagram
\[
\xymatrix{
0\ar[r]&Hom(A,I^{r-1})\ar[d]\ar[r]^-{id\otimes 
i_{*}}&Hom(A,V_{s})\ar[d]^{p^*\otimes id}\ar[r]&Ho
m(A,A)\ar[d]^{\cong}\ar[r]&0\\
0\ar[r]&ad'(V_{s})\ar[r]&Hom(V_{s},V_{s})\ar[r]^-{Tr}&I\ar[r]&0
}
\]
This induces a commutative diagram
\[
\xymatrix{
H^1(Hom(A,I^{r-1}))\ar[r]^{id\otimes i_{*}}\ar[d]
&H^1(Hom(A,V_{s}))\ar[d]^{p^{*}\otimes 
id}\ar[r]&H^1(\mathcal{O}_{Y})\ar[d]^{\cong}\ar[r]&0\\
H^1(ad'(V_{s}))\ar@{^{(}->}[r]&H^1(Hom(V_{s},V_{s}))\ar[r]\ar[d]
&H^1(\mathcal{O}_{Y})\ar[r]&0\\
&0&
}
\]
We conclude that the image of $p^{*}\otimes i_{*}$ equals 
$H^1(ad'(V_{s}))$.
\par
Statement (vi) is immediate from (vii). Indeed, let $f:S\to J$ be the 
vector bundle of (iii). By the construction of \cite{NR1} 
Proposition~2.6 the morphism $f$ is a composition of the  
morphism 
$\det :S\to Pic^{e}Y$ and a translation by a fixed line 
bundle. So, one has a commutative diagram with exact rows
\[
\xymatrix{
0\ar[r]&T_{s}S_{A}\ar[r]\ar[d]&T_{s}S\ar[d]\ar[r]^-{df(s)}
&T_{f(s)}J\ar[d]^{\cong}\ar[r]&0\\
0\ar[r]&H^1(ad'(V_{s}))\ar[r]&H^1(End\, 
V_{s})\ar[r]^-{Tr^*}&H^1(\mathcal{O}_{Y})\ar[r]&0
}
\]
Hence (vii) implies the surjectivity of $\kappa:T_{s}S\to H^1(End\, 
V_{s})$.
\end{proof}
\section{Gorenstein coverings of degree 3,  4 or 5}\label{s5}
In this section we recall and complement some results from \cite{Mi},
\cite{CE} and \cite{C2}. We frequently use Theorem~2.1 of \cite{CE},
whose corrected formulation may be found in \cite{CN} Theorem~2.2. 
In this section, unless otherwise specified, we assume that the base field $k$ is algebraically closed and $char(k) =0$. 
\par
We recall that a complex of sheaves
$\mathcal{F}_{\bullet}: \cdots \to
\mathcal{F}_{i}\overset{d_{i}}
{\lto}\mathcal{F}_{i-1}\to \cdots \to \mathcal{F}_{0}$ is called 
acyclic if $\mathcal{H}_i(\mathcal{F}_{\bullet})=Ker(d_{i})/Im(d_{i+1})=0$ for 
every $i\geq 1$.
\begin{lem}\label{s1.p2}
Let $\rho:P\to Y$ be a flat morphism of schemes. 
\par
A. Let 
\(
\mathcal{F}_{\bullet}: 0\to \mathcal{F}_{s}\to \mathcal{F}_{s-1}\to 
\cdots 
 \to \mathcal{F}_{0}
\)
be a complex of coherent $\mathcal{O}_{P}$-modules, $s\geq 1$. 
 Suppose 
$\mathcal{F}_{0}, \ldots, \mathcal{F}_{s-1}$ are locally free.
Suppose that for every $y\in \rho(P)$ the restriction 
$j^{*}_{y}(\mathcal{F}_{\bullet})$
to the fiber $X_{y}$
is acyclic. Then:
\begin{enumerate}
\item
$\mathcal{F}_{\bullet}$ is acyclic;
\item
$\mathcal{H}_0(\mathcal{F}_{\bullet})=\mathcal{F}_{0}/d_{1}(\mathcal{F}_{1})$ 
is flat over $Y$;
\item
if $j^{*}_{y}(\mathcal{F}_{s})$ is locally free for 
every $y\in \rho(P)$, then $\mathcal{F}_{s}$ is locally free.
\end{enumerate}
\par
B. Let 
\(
 0\to \mathcal{F}_{s}\to \mathcal{F}_{s-1}\to 
\cdots 
 \to \mathcal{F}_{0} \to \mathcal{G} \to 0
\)
be a complex of coherent $\mathcal{O}_{P}$-modules.
Suppose $\mathcal{G}$ is flat over $Y$ and if $s\geq 1$ suppose
$\mathcal{F}_{0}, \ldots, \mathcal{F}_{s-1}$ are locally free. Then the 
complex 
is exact if and only if 
its restriction to each fiber $X_{y}$ is exact.
\end{lem}
\begin{proof}
(A). (i) and (ii) are statements about 
$\mathcal{H}_i(\mathcal{F}_{\bullet})$ which are
proved by descending induction on $s$ using \cite{Ma}~(20.E) or 
\cite{Ma2} Theorem~22.5. They imply that $\mathcal{F}_{s}$ is  flat 
over $Y$. Using \cite{Ma}~(20.G) we conclude, under 
the assumption of (iii), that $\mathcal{F}_{s}$ is a locally free 
$\mathcal{O}_{P}$-module.
\par
(B). The ``if'' part is proved as in (A). For the ``only if'' part one splits the 
exact complex in a sequence of short exact sequences. 
The syzygy sheaves are flat over $Y$ since $\mathcal{G}$ is. 
For  every base change $Y'\to Y$ the inverse image with respect to
the morphism $X\times_{Y}Y'\to X$
preserves the exactness of these short exact sequences because of the flatness
of the right-hand sheaves. Hence it preserves the  exactness of the
complex. The claim follows
considering the base change $\Spec k\to Y$  with image $y\in Y$.
\end{proof}
Let $f:X\to Y$ be a finite, flat, surjective morphism of schemes. Let 
$\omega_{X/Y}$ be its relative dualizing sheaf (see \cite{Kl} or 
\cite{Liu} Ch.6 Proposition~4.25)
\renewcommand{\theenumi}{\roman{enumi}}
\begin{dfn0}
The covering $f:X\to Y$ is called Gorenstein if one of the following 
equivalent conditions is satisfied:
\begin{enumerate}
\item
$\omega_{X/Y}$ is an invertible sheaf;
\item
for every $y\in Y$ the fiber $X_{y}$ is Gorenstein, i.e. 
$\omega_{X_{y}}$ is invertible.
\end{enumerate}
\end{dfn0}
\begin{block}\label{s1.1}
\emph{Triple coverings}\quad \cite{Mi}, \cite{CE}.
\par
Let $Y$ be a scheme.
Let $E$ and $L$ be  locally free sheaves on $Y$ of ranks 2 and 1 
respectively. Consider the 
projective bundle $\rho: \mathbf{P}(E)\to Y$. Let
\[
\Phi_{3}:H^0(Y,S^3E\otimes L)\overset{\sim}{\lto}
H^0(\mathbf{P}(E),\rho^{*}(L)(3))
\]
be the canonical isomorphism. Let $\eta \in 
H^0(Y,S^3E\otimes L)$ and let $\delta=\Phi_{3}(\eta)$.
\begin{dfn0}
The section $\eta$ has  \emph{the right codimension} at $y\in Y$ if 
$\delta|_{\mathbf{P}(E)_{y}}$ has a finite number of zeros
\end{dfn0}
Suppose $\eta \in H^0(Y,S^3E\otimes L)$ has the right 
codimension at every $y\in Y$.
Consider the complex 
\[
0 \lto \rho^{*}(L^{-1})(-3)
\overset{.\delta}{\lto} \mathcal{O}_{\mathbf{P}(E)}.
\]
It is acyclic by Lemma~\ref{s1.p2}, so for every $x\in \mathbf{P}(E)$
the germ $\delta_{x}$ is  
not a zero divisor. Let $X_{\eta}\subset 
\mathbf{P}(E)$ be the Cartier divisor associated with $\delta$. The 
projection $\pi_{\eta}:X_{\eta}\to Y$ is quasifinite and proper, thus 
finite, and surjective. One has an exact complex of 
$\mathcal{O}_{\mathbf{P}(E)}$-modules 
\begin{equation}\label{es1.p4}
0 \lto \rho^{*}(L^{-1})(-3)
\overset{.\delta}{\lto} \mathcal{O}_{\mathbf{P}(E)}\lto 
\mathcal{O}_{X_{\eta}} \lto 0
\end{equation}
and $\mathcal{O}_{X_{\eta}}$ is 
flat over $Y$ by Lemma~\ref{s1.p2}. 
Therefore $\pi_{\eta}:X_{\eta}\to Y$ is a 
covering of degree 3 which is clearly Gorenstein.
\end{block}
\begin{block}\label{s1.3}
\emph{Quadruple coverings}\quad \cite{CE}.
\par
Let $Y$ be a scheme.
Let $E, F$ and $L$ be locally free sheaves on $Y$ of ranks 3,\; 2 and 1 
respectively. 
Consider the projective bundle $\rho 
:\mathbf{P}(E)\to Y$. 
 Let 
\begin{equation*}
\Phi_{4}:H^0(Y,\Check{F}\otimes S^2E \otimes L) \overset{\sim}{\lto} 
H^0(\mathbf{P}(E),\rho^{*}(\Check{F}\otimes L)(2))
\end{equation*}
be the canonical isomorphism. Let $\eta \in H^0(Y,\Check{F}\otimes 
S^2E\otimes L)$ 
and let $\delta = \Phi_{4}(\eta)$.
\begin{dfn0}
The section $\eta$ has  \emph{the right codimension} at $y\in Y$ if 
the pencil 
of conics in $\mathbf{P}(E)_{y}$, determined by $\eta(y):(F\otimes 
L^{-1})\otimes k(y)\to (S^2E)\otimes k(y)$, has base locus of dimension 0. 
\end{dfn0}
Suppose $\eta \in 
H^0(Y,\Check{F}\otimes S^2E \otimes L)$ has the right codimension at every 
$y\in Y$. Consider the Koszul complex associated with the section 
$\delta \in H^0(\mathbf{P}(E),\rho^{*}(\Check{F}\otimes L)(2))$ 
(cf. \cite{FL}~pp.75,76)
\begin{equation}\label{es1.3aa}
0 \lto \rho^{*}\bigwedge^{2}(F\otimes L^{-1})(-4) \lto \rho^{*}(F\otimes 
L^{-1})(-2) \lto \mathcal{O}_{\mathbf{P}(E)}.
\end{equation}
Let $y\in Y$. Restricting on $\mathbf{P}(E)_{y}$ one 
obtains an acyclic complex since $\mathbf{P}(E)_{y}$ is nonsingular, 
thus Cohen-Macaulay,  and the zero locus of $\delta 
|_{\mathbf{P}(E)_{y}}$ has codimension 2. By Lemma~\ref{s1.p2} the 
complex \eqref{es1.3aa} is acyclic. Let $X_{\eta}\subset \mathbf{P}(E)$ be the 
zero subscheme of $\delta$. As in the case $d=3$ the 
projection $\pi_{\eta}:X_{\eta}\to Y$ is finite, flat and surjective, and one 
has the 
following resolution of $\mathcal{O}_{X_{\eta}}$ by locally free 
$\mathcal{O}_{\mathbf{P}(E)}$-modules:
\begin{equation} \label{es1.p5}
0 \lto \rho^{*}(\det F\otimes L^{-2})(-4) \lto \rho^{*}(F\otimes 
L^{-1})(-2) \lto \mathcal{O}_{\mathbf{P}(E)} \lto 
\mathcal{O}_{X_{\eta}} \lto 0.
\end{equation}
The covering $\pi_{\eta}:X_{\eta}\to Y$ is clearly Gorenstein of 
degree 4 and furthermore the embedding $i:X_{\eta}\hookrightarrow 
\mathbf{P}(E)$ has the property that every fiber 
$(X_{\eta})_{y}\subset\mathbf{P}(E)_{y}$, being a zerodimensional complete 
intersection of two conics, is arithmetically Gorenstein and nondegenerate. 
 
\end{block}
We use below some results of Buchsbaum and Eisenbud \cite{BuE}. We work over 
a field of characteristic zero, so the use of divided powers may be avoided. For reader's
convenience we include simplified proofs of two statements due to Buchsbaum and Eisenbud.
We recall some facts about exterior algebras and the operations $\wedge$ and $\lrcorner$
for which we refer to \cite{St} Chapter~1.
Let $V$ be a vector space over a field $k$ of characteristic zero.  There is a canonical
duality $\langle \; | \; \rangle : \bigwedge ^{\bullet}V^*\times \bigwedge ^{\bullet}V\to k$ such that 
$\langle a^1\wedge \cdots \wedge a^p|b_1\wedge \cdots \wedge b_p \rangle = \det(\langle a^i , b_j \rangle )$, where $a^i\in V^*, b_j\in V$. If $u^*\in \bigwedge^pV^*$ and $x\in \bigwedge^nV$, then $u^*\lrcorner x\in 
\bigwedge^{n-p}V$ and it is defined by $\langle v^*|u^*\lrcorner x\rangle = 
\langle u^* \wedge v^*| x\rangle$. If $R$ is a commutative $k$-algebra and $F=R\otimes_kV$ the above formulas extend to a duality $\bigwedge ^{\bullet}F^*\times \bigwedge ^{\bullet}F\to R$ and $R$-bilinear operations $\wedge$ and $\lrcorner$. If $\varphi \in \bigwedge^2F$, then $\wedge^m\varphi := \varphi \wedge \cdots \wedge \varphi \in \bigwedge^{2m}F$.
\begin{lem}[Buchsbaum-Eisenbud]\label{BE2}
Let $R$ be a commutative $k$-algebra, where $k$ is a field of characteristic zero. Let $F$ be a free $R$-module of odd rank $n=2m+1\geq 3$. Let $e_1,\ldots,e_n$ be its basis and let $e^1,\ldots,e^n$ be the dual basis of $F^*$. Let $\varphi = \sum_{1\leq i<j\leq n}x_{ij}e_i\wedge e_j \in \bigwedge^2F$. Then one has the following complex
\begin{equation} \label{eBE2}
0 \to \bigwedge^nF^* 
\overset{^{t}g}{\lto} F^* \overset{f}{\lto} F
\overset{g}{\lto} \bigwedge^nF
\end{equation}
where $f(a^*)= a^*\lrcorner \varphi$, $g(b)=b\wedge (\wedge^m\varphi)$, 
$^{t}g(e^1\wedge \cdots \wedge e^n)=\wedge^m\varphi \lrcorner e^1\wedge \cdots \wedge e^n$. Let $\Phi = (x_{ij})_{i,j=1}^n$, where $x_{ji}=-x_{ij}$. Then the matrix of $g$ is
\[
m!(Pf(\Phi_1),\ldots,(-1)^{i-1}Pf(\Phi_i),\ldots,(-1)^{n}Pf(\Phi_n)),
\]
where $Pf(\Phi_i), i=1,\ldots,n$ are the submaximal pfaffians, i.e. $\Phi_i$ is the skew-symmetric matrix obtained by deleting the $i$-th row and the $i$-th column of $\Phi$ and $Pf(\Phi_i)$ is its pfaffian.
\end{lem}
\begin{proof}
$D_{a^*}(u)=a^*\lrcorner u$ is an antiderivation of $\bigwedge^{\bullet}F$ (see \cite{St} Ch.~1 Ex.~4.5). Therefore $0=a^*\lrcorner \wedge^{m+1}\varphi = (m+1)(a^*\lrcorner \varphi)\wedge (\wedge^{m}\varphi)$. This proves $g\circ f=0$. Transposing, one obtains $^{t}f\circ {^{t}g}=0$. One has $^{t}f=-f$ and $^{t}g:(\bigwedge^nF)^*\to F^*$ is given by
\(
\langle ^{t}g(e^1\wedge \cdots \wedge e^n)| b\rangle = \langle e^1\wedge \cdots \wedge e^n| b\wedge (\wedge^m\varphi)\rangle = \langle \wedge^m\varphi \lrcorner e^1\wedge \cdots \wedge e^n| b\rangle
\).
This shows that $f\circ {^{t}g} = 0$, so  \eqref{eBE2} is a complex.
\par
One has, according to  \cite{Bou} \S~5 n.2, that 
\begin{equation*}
\begin{split}
g(e_i) &=e_i\wedge (\wedge^m\varphi) 
       = e_i\wedge (\wedge^m\sum_{s<t, i\notin \{s,t\}}x_{st}e_s\wedge e_t) \\
       &= m!e_i\wedge Pf(\Phi_i)e_1\wedge\cdots\wedge e_{i-1}\wedge e_{i+1}\wedge\cdots\wedge e_n \\
       &= m!(-1)^{i-1} Pf(\Phi_i)e_1\wedge\cdots\wedge e_n.
\end{split}
\end{equation*}
This shows the last statement of the lemma.
\end{proof}
\begin{pro}[Buchsbaum-Eisenbud]\label{BE3}
Let $X$ be a connected Cohen-Macaulay scheme of dimension $\geq 3$. Let $F$ be a locally free sheaf on $X$ of rank 
$n=2m+1\geq 3$ and let $M$ be an invertible sheaf on $X$. Let $\varphi \in H^0(X,(\bigwedge^2F) \otimes M)$. Let $Y\subset X$ be the closed subscheme whose ideal sheaf is generated by the submaximal pfaffians associated with $\varphi$. Suppose $\{x\in X|\rank \varphi(x) < 2m\}$, which as a set coincides with $Y$, has codimension $\geq 3$ in $X$. Then $\mathcal{O}_Y$ has the following resolution by locally free sheaves:
\begin{equation}\label{eBE4a}
\begin{aligned}
0&\lto M^{-n}\otimes (\bigwedge^nF)^{-2}\lto F^{\vee}\otimes M^{-m-1}\otimes 
(\bigwedge^nF)^{-1} \\
 &\lto F\otimes M^{-m}\otimes (\bigwedge^nF)^{-1}\lto \mathcal{O}_X\lto \mathcal{O}_Y\lto 0.
\end{aligned}
\end{equation}
Furthermore $Y$ is an equidimensional Cohen-Macaulay subscheme of codimension 3.
\end{pro}
\begin{proof}
Using Lemma~\ref{BE2} one obtains the following complex:
\begin{equation}\label{eBE4b}
0\lto \bigwedge^nF^{\vee}\otimes (M^{\vee})^{m}\overset{^{t}g}{\lto} F^{\vee}  
\overset{f}{\lto} F\otimes M \overset{g}{\lto} \bigwedge^nF \otimes M^{m+1}
\end{equation}
where the differentials are defined on local sections as: 
\begin{equation*}
\begin{aligned}
&f(a^*)= a^*\lrcorner \varphi,\quad g(b)=b\wedge (\wedge^m\varphi), \\
&{^{t}g}(f^1\wedge \cdots \wedge f^n\otimes (e^*)^{\otimes (m)})=\wedge^m\varphi \lrcorner (f^1\wedge \cdots \wedge f^n \otimes (e^*)^{\otimes (m)}).
\end{aligned}
\end{equation*}
We claim the complex \eqref{eBE4b} is acyclic. This is a local matter, so we may assume that: $X$ is an affine scheme;  $F$ and $M$ are free sheaves. Let $X = \Spec R$ and let $I$ be the ideal generated by the submaximal pfaffians of the skew-symmetric matrix associated with $\varphi$. According to the Buchsbaum-Eisenbud acyclicity criterion \cite{BE2}, which we apply in the form given in \cite{BrH} Theorem~1.4.13, one has to verify that: a) $\grade I_1({^{t}g})\geq 3$; b) $\grade I_{2m}(f)\geq 2$; c) $\grade I_1(g)\geq 1$ (we notice that $\grade J=\depth(J,R)$ in the terminology of \cite{Ei}). By Lemma~\ref{BE2} one has $I = I_1(g)=I_1({^{t}g})$.  The three ideals $I_1({^{t}g}), I_{2m}(f)$ and $I_1(g)$ determine one and the same closed  subset $\{x\in X|\rank \varphi(x) < 2m\}$ since for a bivector $\omega = \frac{1}{2}\sum_{i,j}a_{ij}e_i\wedge e_j \in \bigwedge^2k^{2m+1}$, where $a_{ji}=-a_{ij}$, the following conditions are equivalent: (i) $\rank (a_{ij})<2m$; (ii) the rank of $\omega$ is less than $2m$; (iii) $\wedge^m\omega = 0$; (iv) the submaximal pfaffians of $(a_{ij})$ are equal to zero. Therefore the radicals of $I_1({^{t}g}), I_{2m}(f)$ and $I_1(g)$ coincide and are equal to $I(Y_{red})\subset R$. The grade of an ideal equals the grade of its radical \cite{Ei} Corollary~17.8, so all three ideals have the same grade, which equals $\height I(Y_{red})= \codim Y$, since $X$ is Cohen-Macaulay. The hypothesis $\codim Y\geq 3$ imples that conditions (a), (b) and (c) above hold. Therefore the complex \eqref{eBE4b} is acyclic. Tensoring \eqref{eBE4b} by the dual of $\bigwedge^nF\otimes M^{m+1}$ one obtains the resolution \eqref{eBE4a}.
\par
For the last statement of the proposition we may again assume $X$ is an affine scheme and $F$ and $M$ are free sheaves. Let $R$ and $I$ be as above. The scheme $X$ is equidimensional being connected and Cohen-Macaulay \cite{Ei} Corollary~18.11. The section $\varphi \in H^0(X,\bigwedge^2\mathcal{O}_X^{\oplus(2m+1)})$ determines a morphism $h:X\to \bigwedge^2k^{2m+1}$ and as a set of points $Y=h^{-1}(\Delta_{2m-2})$, where $\Delta_{2m-2}\subset \bigwedge^2k^{2m+1}$ is the locus of bivectors of rank $\leq 2m-2$. Since $\codim \Delta_{2m-2} = 3$ we conclude that at every point $y\in Y$, $\codim_yY\leq 3$. Hence the hypothesis $\codim Y\geq 3$ implies $Y$  is equidimensional of codimension 3. The exact complex \eqref{eBE4a} and \cite{Ei} Corollary~18.5 show that $3\geq \projdim_{R}R/I\geq \grade I = 3$. Hence $I$ is a perfect ideal and therefore $R/I$ and $Y=\Spec R/I$ are Cohen-Macaulay (see \cite{BrH} Theorem~2.1.5).
\end{proof}
\begin{block}\label{s5.80}
\emph{Coverings of degree 5}\quad \cite{BuE},\: \cite{C2}.
\par
Let $Y$ be a scheme. 
Let $E$, $F$ and $L$ be locally free 
sheaves on $Y$ of ranks 4, 5 and 1 respectively. 
 Consider the projective bundle $\rho: \mathbf{P}(E)\to 
Y$ and the canonical isomorphism 
\[
\Phi_{5}:H^0(Y,\bigwedge^2F\otimes E\otimes L)
\overset{\sim}{\lto} H^0(\mathbf{P}(E),\rho^*((\bigwedge^2F)\otimes L)
(1))
\]
Let $\eta \in H^0(Y,\bigwedge^2F\otimes 
E\otimes L)$ and let $\delta = \Phi_{5}(\eta)$. For every 
$y\in Y$ the restriction $\delta|_{\mathbf{P}(E)_{y}}$ belongs to 
$H^0(\mathbb{P}^{3},\bigwedge^{2}\mathbb{C}^{5}\otimes 
\mathcal{O}_{\mathbb{P}^{3}}(1))$, where 
$\mathbb{P}^{3}=\mathbf{P}(E)_{y}$. So, $\delta|_{\mathbf{P}(E)_{y}}$ is a
bivector,
 whose coefficients are linear forms of $\mathbb{P}^{3}$. The 
locus of points in $\mathbb{P}^{3}$ where this bivector has 
rank $\leq 2$, is closed of codimension $\leq 3$ in $\mathbf{P}(E)_{y}$. 
\begin{dfn0}
The section  $\eta \in 
H^0(Y,\bigwedge^2F\otimes E\otimes L)$
 has \emph{the right codimension} at 
$y\in Y$, if $\delta|_{\mathbf{P}(E)_{y}}$ has rank $\leq 2$ only at a 
finite number of 
points in $\mathbf{P}(E)_{y}$.
\end{dfn0}
\begin{rem0}
We notice that this condition is weaker than the one in 
\cite{C2} Definition~3.6. We do not require that $\delta$ vanishes in 
no point of $\mathbf{P}(E)_{y}$. 
\end{rem0}
Let $\eta \in H^0(Y,\bigwedge^2F\otimes E\otimes L)$ and let $\delta = 
\Phi_{5}(\eta)$. If $U\subset Y$ is an open subset such that $E|_{U}, 
F|_{U}$ and $L|_{U}$ are trivial, then $\delta|_{\rho^{-1}(U)}$ determines
 a skew-symmetric $5\times 5$ matrix with entries in 
$\Gamma (\rho^{-1}(U),\mathcal{O}_{\mathbf{P}(E)}(1))$. Let 
$\phi_{1},\ldots,\phi_{5}$ be the pfaffians of the five $4\times 4$ 
principal minors of the matrix. They belong to 
$\Gamma (\rho^{-1}(U),\mathcal{O}_{\mathbf{P}(E)}(2))$. So, 
$\sum_{i=1}^{5} \mathcal{O}_{\rho^{-1}(U)}(-2)\cdot \phi_{i}$ is an ideal 
subsheaf of $\mathcal{O}_{\rho^{-1}(U)}$. Taking an open covering of 
$Y$ and patching, one obtains an ideal subsheaf of 
$\mathcal{O}_{\mathbf{P}(E)}$. It defines a closed subscheme 
$X_{\eta}\subset \mathbf{P}(E)$. Let $\pi_{\eta}:X_{\eta}\to Y$ be the 
restriction $\pi_{\eta} = \rho|_{X_{\eta}}$. 
\end{block}
\begin{pro}\label{s5.81}
Let the assumptions be as in \S\ref{s5.80}.
\par 
\quad (i) Suppose 
$\eta \in 
H^0(Y,\bigwedge^2F\otimes E\otimes L)$
 has the right codimension at every $y\in Y$. Let 
$X=X_{\eta}\subset\mathbf{P}(E)$. Then $\pi = \rho|_{X}:X\to Y$ is a 
Gorenstein covering of degree 5 such that every fiber $X_{y}\subset 
\mathbf{P}(E)_{y}$ is non degenerate and arithmetically Gorenstein. 
The structure sheaf $\mathcal{O}_{X}$ has the following resolution:
\begin{equation}\label{es5.81}
0 \to \rho^{*}(F_{3})(-5) \to \rho^{*}(F_{2})(-3) \to 
\rho^{*}(F_{1})(-2) \to \mathcal{O}_{\mathbf{P}(E)} \to 
\mathcal{O}_{X} \to 0
\end{equation}
where $F_{3}=(\det F)^{-2}\otimes L^{-5}$, $F_{2}=\Check{F}\otimes 
(\det F)^{-1}\otimes L^{-3}$, $F_{1}=F\otimes (\det F)^{-1}\otimes 
L^{-2}$.
 \par  \medskip
\quad (ii) Suppose furthermore that $Y$ is integral and
$L\cong (\det E)^{-1}$, $\det F 
\cong (\det E)^{2}$. Then $E\cong 
(\pi_{*}\mathcal{O}_{X}/\mathcal{O}_{Y})^{\vee}$ and the 
closed embedding of \,
$Y$-schemes $X\hookrightarrow \mathbf{P}(E)$ is the composition 
\(
\xymatrix{
X\ar[r]&\mathbf{P}(\pi_{*}\omega_{X/Y})\ar@{-->}[r]&\mathbf{P}(E)
}
\)
induced by the exact sequence
\[
0 \lto (\pi_{*}\mathcal{O}_{X}/\mathcal{O}_{Y})^{\vee} \lto 
\pi_{*}\omega_{X/Y}\overset{t_{\pi}}{\lto} \mathcal{O}_{Y}\lto 0
\]
where $t_{\pi}:\pi_{*}\omega_{X/Y}\to \mathcal{O}_{Y}$ is the trace map. 
Moreover the resolution \eqref{es5.81} of (i) becomes
\begin{equation}\label{es5.81a}
0 \to \rho^{*}(\det E)(-5) \to \rho^{*}(\Check{F}\otimes \det E)(-3) \to 
\rho^{*}(F)(-2) \to \mathcal{O}_{\mathbf{P}(E)} \to 
\mathcal{O}_{X} \to 0
\end{equation}
and is isomorphic to the canonical resolution of $\mathcal{O}_{X}$ 
defined in \cite{CE} Theorem~2.1.
\par \medskip
\quad (iii) Suppose $Y$ is integral. Given the data of (i) 
one may replace $(\eta,F,E,L)$ by 
a quadruple $(\lambda\eta,F',E',L')$ which satisfies 
the conditions of (ii). 
 Here $F'=F\otimes T_{1}$, $E'=E\otimes 
T_{2}$, $L'=L\otimes T_{3}$, where $T_{1},T_{2},T_{3}$ are invertible 
sheaves on $Y$ which
satisfy
$T_{1}^{2}T_{2}T_{3}\cong \mathcal{O}_{Y}$, and $\lambda$ is a 
trivializing section of $T_{1}^{2}T_{2}T_{3}$. This normalization 
preserves $X\subset \mathbf{P}(E)$. The triple $T_{1},T_{2},T_{3}$
is unique up to isomorphism.
\end{pro}
\begin{proof}
(i) Let $P_{2} = \rho^{*}(\Check{F})$, $P_{3}=\rho^{*}(L)(1)$. We may 
identify $\delta$ with  $\varphi \in H^0(\mathbf{P}(E),
\bigwedge^{2}P_{2}^{\vee}\otimes P_{3})$. One has the Buchsbaum-Eisenbud complex 
\eqref{eBE4b}
\begin{equation}\label{es5.82}
0 \to \bigwedge^{5}P_{2}\otimes (P_{3}^{\vee})^{\otimes 2} 
\overset{^{t}g}{\lto} P_{2} \overset{f}{\lto} P_{2}^{\vee}\otimes P_{3}
\overset{g}{\lto} \bigwedge^{5}P_{2}^{\vee}\otimes P_{3}^{\otimes 3}
\end{equation}
where the differentials are defined on local sections as
\begin{equation}\label{es5.82a}
\begin{aligned}
&f(a) = a \lrcorner \varphi , \qquad
g(b) =b\wedge (\wedge^{2}\varphi) \\
&^{t}g(f_{1}\wedge \cdots \wedge f_{5}\otimes (e^{*})^{\otimes 2}) =
\wedge^{2}\varphi \lrcorner 
(f_{1}\wedge \cdots \wedge f_{5}\otimes (e^{*})^{\otimes 2}).
\end{aligned}
\end{equation}
Let $y\in Y$. Tensoring \eqref{es5.82} by 
$\mathcal{O}_{\mathbf{P}(E)_{y}}$ one obtains an acyclic complex. 
Here one applies Proposition~\ref{BE3} taking into account that $\mathbf{P}(E)_{y}$ is nonsingular, thus Cohen-Macaulay. By Lemma~\ref{s1.p2} the complex \eqref{es5.82}
 is acyclic and $Coker(g)$ is flat over $Y$. Tensoring by 
$\bigwedge^{5}P_{2}\otimes (P_{2}^{\vee})^{\otimes 3}$ one obtains the 
resolution \eqref{es5.81}. As in the cases $d=3$ and $d=4$ the 
projection morphism $\pi_{\eta}:X_{\eta}\to Y$ is finite, flat and surjective. 
The fibers are Gorenstein schemes, being complete intersections
of a Fano variety - the Pl\"{u}cker embedding of a Grassmanian.
Therefore $\pi_{\eta}:X_{\eta}\to Y$ is a
Gorenstein covering of degree 5. The embedding 
$X_{\eta}\hookrightarrow \mathbf{P}(E)$ has the property that 
every fiber $(X_{\eta})_{y}\subset \mathbf{P}(E)_{y}$ is arithmetically 
Gorenstein and nondegenerate. The latter is true since the 
Buchsbaum-Eisenbud complex yields a graded minimal resolution of the 
coordinate ring $S((X_{\eta})_{y})$, so the homogeneous ideal of 
$(X_{\eta})_{y}$ is generated by forms of degree 2. 
\par
(ii) Here the resolution \eqref{es5.81} equals \eqref{es5.81a}. In 
particular the last term is isomorphic to $\rho^{*}(\det E)(-5)$. One 
applies Part~(iv) of Theorem~2.1 of \cite{CE}. 
\par
(iii) Let the invertible sheaves $T_1,T_2,T_3$ satisfy 
$T_1^{2}T_2T_3\cong \mathcal{O}_{Y}$ and let $\lambda \in 
H^0(Y,T_1^2T_2T_3)$ be a trivializing section. Replacing 
$(\eta,F,E,L)$ by $(\lambda\eta,F\otimes T_{1},E\otimes T_{2},
L\otimes T_{3})$ one has $\mathbf{P}(E)=\mathbf{P}(E\otimes T_{2})$ 
and $X_{\eta} = X_{\lambda \eta}$.
\par
As shown in Part~C of the proof of Theorem~2.1 of \cite{CE} there is a 
unique $[M]\in Pic\, Y$ such that replacing $E$ by $E'=E\otimes M$, 
$\mathcal{O}_{\mathbf{P}(E)}(1)$ by $\mathcal{O}_{\mathbf{P}(E')}(1) = 
\mathcal{O}_{\mathbf{P}(E)}(1)\otimes \rho^{*}M$ and correspondingly 
modifying the resolution \eqref{es5.81} one obtains for the last term 
$F'_{3}\cong \det E'$. Furthermore $M$ is isomorphic to $\Check{F}_{3}
\otimes \det E$. We let $T_{2}=(\det F)^{2}\otimes  L^{ 5}\otimes \det E$. 
One  imposes one additional condition  on  $T_1,T_2,T_3$ 
in order to have $(\det F')\otimes  L'^{2}  \cong  \mathcal{O}_{Y}$,  so 
that $F'_{1}\cong F'$. One obtains the following system  of  equations 
in $Pic\, Y$:
\begin{equation}\label{es5.83}
\begin{split}
T_1^2T_2T_3 &\cong \mathcal{O}_{Y}\\
T_{2}       &\cong (\det F)^{2}\otimes L^{5}\otimes 
\det E\\
T_{1}^{5}T_{3}^{2}   &\cong   (\det   F)^{-1}\otimes 
L^{-2}. 
\end{split}
\end{equation}
The system has obviously a unique solution up to isomorphisms of
$T_i, i=1,2,3$. Consider the resolution of 
$\mathcal{O}_{X}$ obtained replacing 
$(\eta,F,E,L)$ by $(\lambda\eta,F',E',L')$. Then $\det E' \cong 
F'_{3}  \cong  (\det  F')^{-2}\otimes  L'^{-5}$  and  $\det  F'\otimes 
L'^{2}\cong \mathcal{O}_{Y}$. Therefore  $\det  E'\cong  L'^{-1}$  and 
$\det F' \cong (\det E')^{2}$. 
\par
The uniqueness statement is true since the system of equations 
\[
t_{1}^{5}=t_{2}^{8},\qquad t_{3}=t_{2}^{-4}, \qquad 
t_{1}^{2}t_{2}t_{3} = 1
\]
has a unique solution $t_{1}=t_{2}=t_{3}=1$ in the group $Pic\, Y$. 
Indeed, eliminating $t_{3}$ one obtains $t_{1}^{5}=t_{2}^{8}, 
t_{1}^{2}=t_{2}^{3}$. Hence $t_{1}=t_{2}^{2} \Rightarrow t_{2}=1 
\Rightarrow t_{1}=1 \Rightarrow t_{3}=1$.
\end{proof}
With similar and simpler arguments, left to the reader, one may 
complement \S~\ref{s1.1} and \S~\ref{s1.3}.
\begin{pro}\label{s5.83a}
Let the assumptions be as in \S~\ref{s1.1}.
\par
(i) Suppose $\eta \in H^0(Y,S^3E\otimes L)$ has the right codimension 
at every $y\in Y$. Let $X=X_{\eta}\subset \mathbf{P}(E)$. Then $\pi = 
\rho|_{X} : X\to Y$ is a Gorenstein covering of degree 3. The 
resolution of $\mathcal{O}_{X}$ is given by \eqref{es1.p4}.
\par
(ii) Suppose furthermore that $Y$ is integral and $L \cong  (\det E)^{-1}$. 
Then 
$E \cong (\pi_{*}\mathcal{O}_{X}/\mathcal{O}_{Y})^{\vee}$ and 
the closed embedding of \,
$Y$-schemes $X\hookrightarrow \mathbf{P}(E)$ is the composition 
\(
\xymatrix{
X\ar[r]&\mathbf{P}(\pi_{*}\omega_{X/Y})\ar@{-->}[r]&\mathbf{P}(E)
}
\)
induced by the exact sequence
\[
0 \lto (\pi_{*}\mathcal{O}_{X}/\mathcal{O}_{Y})^{\vee} \lto 
\pi_{*}\omega_{X/Y}\overset{t_{\pi}}{\lto} \mathcal{O}_{Y}\lto 0
\]
where $t_{\pi}:\pi_{*}\omega_{X/Y}\to \mathcal{O}_{Y}$ is the trace map. 
The resolution of $\mathcal{O}_{X}$ equals
\begin{equation*}
0 \to \rho^{*}(\det E)(-3) \to \mathcal{O}_{\mathbf{P}(E)} \to 
\mathcal{O}_{X} \to 0
\end{equation*}
and is isomorphic to the canonical resolution of $\mathcal{O}_{X}$ 
defined in \cite{CE} Theorem~2.1.
\par
(iii) Suppose $Y$ is integral. Given the data of (i) one may replace $(\eta, 
E, 
L)$ by a triple 
$(\lambda\eta, E', L')$ which satisfies the conditions of (ii). 
  Here  $E'=E\otimes 
T_{1}$, $L'=L\otimes T_{2}$, where $T_{1},T_{2}$ are invertible 
sheaves on $Y$ which
satisfy
$T_{1}^{3}T_{2}\cong \mathcal{O}_{Y}$, and $\lambda$ is a 
trivializing section of $T_{1}^{3}T_{2}$. This normalization 
preserves $X\subset \mathbf{P}(E)$. The couple $T_{1},T_{2}$
is unique up to isomorphism.
\end{pro}
\begin{proof}
We only mention that the analog of \eqref{es5.83} is: $T_1^3T_2\cong \mathcal{O}_Y; T_1\cong L\otimes \det E$.
\end{proof}
\begin{pro}\label{s5.83b}
Let the assumptions be as in \S~\ref{s1.3}.
\par
(i) Suppose $\eta \in H^0(Y,\Check{F}\otimes S^2E\otimes L)$ has the 
right codimension at every $y\in Y$. Let $X=X_{\eta}\subset 
\mathbf{P}(E)$. Then $\pi=\rho|_{X}:X\to Y$ is a Gorenstein covering 
of degree 4. Furthermore every fiber $X_{y}\subset \mathbf{P}(E)_{y}$ 
is nondegenerate and arithmetically Gorenstein. The resolution of 
$\mathcal{O}_{X}$ is given by \eqref{es1.p5}. 
\par
(ii) Suppose furthermore that $Y$ is integral and $L\cong \mathcal{O}_{Y}$, 
$\det F\cong 
\det E$. 
Then 
$E \cong (\pi_{*}\mathcal{O}_{X}/\mathcal{O}_{Y})^{\vee}$ and 
the closed embedding of \,
$Y$-schemes $X\hookrightarrow \mathbf{P}(E)$ is the composition 
\(
\xymatrix{
X\ar[r]&\mathbf{P}(\pi_{*}\omega_{X/Y})\ar@{-->}[r]&\mathbf{P}(E)
}
\)
induced by the exact sequence
\[
0 \lto (\pi_{*}\mathcal{O}_{X}/\mathcal{O}_{Y})^{\vee} \lto 
\pi_{*}\omega_{X/Y}\overset{t_{\pi}}{\lto} \mathcal{O}_{Y}\lto 0
\]
where $t_{\pi}:\pi_{*}\omega_{X/Y}\to \mathcal{O}_{Y}$ is the trace map. 
The resolution of $\mathcal{O}_{X}$ equals
\begin{equation*}
0 \lto \rho^{*}(\det E)(-4) \lto \rho^{*}(F)(-2) \lto 
\mathcal{O}_{\mathbf{P}(E)} \lto 
\mathcal{O}_{X_{\eta}} \lto 0.
\end{equation*}
and is isomorphic to the canonical resolution of $\mathcal{O}_{X}$ 
defined in \cite{CE} Theorem~2.1.
\par
(iii) Given the data of (i) one may replace $(\eta,F,E,L)$ by 
a quadruple $(\lambda\eta,F',E',L')$ which satisfies 
the conditions of (ii). 
 Here $F'=F\otimes T_{1}$, $E'=E\otimes 
T_{2}$, $L'=L\otimes T_{3}$, where $T_{1},T_{2},T_{3}$ are invertible 
sheaves on $Y$ which
satisfy
$T_{1}^{-1}T_{2}^{2}T_{3}\cong \mathcal{O}_{Y}$, and $\lambda$ is a 
trivializing section of $T_{1}^{-1}T_{2}^{2}T_{3}$. This normalization 
preserves $X\subset \mathbf{P}(E)$. The triple $T_{1},T_{2},T_{3}$
is unique up to isomorphism.
\end{pro}
\begin{proof}
We only mention that here the analog of \eqref{es5.83} is
\begin{equation*}
\begin{split}
T_1^{-1}T_2^{2}T_3 &\cong \mathcal{O}_{Y}\\
T_{2}       &\cong (\det F)^{-1}\otimes \det E\otimes L^{2}
\\
T_{3}   &\cong   L^{-1}. 
\end{split}
\end{equation*}
\end{proof}
\begin{cor}\label{s5.83c}
Let $Y$ be an irreducible variety. Let $d=3,4$ or $5$. Let $\eta$ be 
a global section of $S^3E\otimes L$, $\Check{F}\otimes S^2E\otimes L$ 
or $\bigwedge^2F\otimes E\otimes L$ respectively. Suppose $\eta$ has 
the right codimension at every $y\in Y$. Let $X=X_{\eta}\subset 
\mathbf{P}(E)$ and let $\pi:X\to Y$ be the associated 
Gorenstein covering of degree $d$. Then 
\begin{equation*}
\omega_{X/Y} \cong
\begin{cases}
\mathcal{O}_{X}(1)\otimes \pi^{*}(L\otimes \det E)          & \text{if 
$d=3$}\\
\mathcal{O}_{X}(1)\otimes \pi^{*}(L^{2}\otimes (\det F)^{-1}\otimes \det E)    
        & \text{if $d=4$}\\
\mathcal{O}_{X}(1)\otimes \pi^{*}(L^{5}\otimes (\det F)^{2}\otimes \det E)
        & \text{if $d=5$}.
\end{cases}
\end{equation*}
\end{cor}
\begin{proof}
In each of the cases the resolution of $\mathcal{O}_{X}$ is given by 
\eqref{es1.p4}, \eqref{es1.p5} and \eqref{es5.81} respectively. The 
last term of this resolution is $\rho^{*}(F_{d-2})(-d)$. Let 
$M=\Check{F}_{d-2}\otimes \det E$. It is proved in \cite{CE} p.447
that, letting $E'=E\otimes M$, the invertible sheaf 
$\mathcal{O}_{\mathbf{P}(E')}(1)\cong \mathcal{O}_{\mathbf{P}(E)}(1)
\otimes \rho^{*}M$ restricts to $\omega_{X/Y}$. One has $F_{d-2}= 
L^{-1}$, $\det F\otimes L^{-2}$ and $(\det F)^{-2}\otimes L^{-5}$ 
when $d=3, 4$ or $5$ respectively. This gives the formula for 
$\omega_{X/Y}$.
\end{proof}
We need the openness of some conditions for families of coverings
analogous to \cite{Gr} Th\'{e}or\`{e}me~12.2.4.
\renewcommand{\theenumi}{\alph{enumi}}
\begin{pro}\label{s5.69}
Let $\mathcal{X}\overset{\pi}{\lto}\mathcal{Y}\overset{q}{\lto}S$ be 
morphisms of schemes, where $q$ is proper, flat, with equidimensional 
fibers of  fixed dimension $n\geq 1$ and $\pi$ is  finite, flat and 
surjective
 of degree $d\geq 2$. The following conditions on $s\in S$ determine 
open (possibly empty) subsets of $S$ (in conditions (g) and (h) below 
it is assumed that
  $\mathcal{X}_{s}\to \mathcal{Y}_{s}$ is 
branched for every $s\in S$):
\begin{enumerate}
\item
$\mathcal{Y}_{s}$ is smooth;
\item
$\mathcal{Y}_{s}$ is irreducible and generically reduced;
\item
$\mathcal{X}_{s}$ is smooth;
\item
$\mathcal{X}_{s}$ is irreducible and generically reduced;
\item
$\pi_{s}:\mathcal{X}_{s}\to \mathcal{Y}_{s}$ is unramified over a 
dense, open subset of $\mathcal{Y}_{s}$;
\item
$\pi_{s}:\mathcal{X}_{s}\to \mathcal{Y}_{s}$ is \'{e}tale;
\item
$\mathcal{Y}_{s}$ has no embedded components and the discriminant scheme of \linebreak
$\pi_{s}:\mathcal{X}_{s}\to 
\mathcal{Y}_{s}$ is of pure codimension one and  smooth;
\item
$\mathcal{Y}_{s}$ has no embedded components and the discriminant scheme of \linebreak
 $\pi_{s}:\mathcal{X}_{s}\to 
\mathcal{Y}_{s}$ is of codimension one, irreducible and generically 
reduced.
\end{enumerate}
\end{pro}
\begin{proof}
(a) $y\in \mathcal{Y}$ is a smooth point of the fiber 
$\mathcal{Y}_{s}$, $s=q(y)$ if and only if $q$ is smooth at $y$ 
(\cite{AK} Theorem~VII.1.8). Moreover such points form an open subset 
of $\mathcal{Y}$ (\cite{AK} Remark~VII.1.2). Let $T$ be its 
complement. Then $q(T)$ is closed in $S$ since $q$ is proper. The open 
set $S\setminus q(T)$ consists of the points which satisfy 
(a).
\par
(b) This is proved in 
\cite{Be} Proposition~1.1.
\par
(c) and (d) The composition $g=q\circ \pi:\mathcal{X}\to S$ is proper, 
 flat, with equidimensional fibers. One repeats the arguments of (a) and (b).
\par
(e) Let $d_{\mathcal{X}/\mathcal{Y}}:(\bigwedge^{max}\pi_{*}
\mathcal{O}_{\mathcal{X}})^{\otimes 2}\lto \mathcal{O}_{\mathcal{Y}}$ be the 
discriminant morphism (cf. \cite{AK} p.124) and let $D_{\mathcal{X}/\mathcal{Y}}$
be its image, the discriminant ideal sheaf. Let us denote by $B=V(D_{\mathcal{X}/\mathcal{Y}})$
the discriminant locus and by $\mathcal{B}$ the discriminant scheme, which is the topological space $B$
endowed by the structure sheaf $\mathcal{O}_{\mathcal{Y}}/D_{\mathcal{X}/\mathcal{Y}}$. Since $D_{\mathcal{X}/\mathcal{Y}}$
is a locally principal ideal sheaf, by Krull's principal ideal theorem, the following possibilities may occur
for a fiber of $q:\mathcal{Y}\to S$:
(i) $B\cap \mathcal{Y}_{s}=\emptyset$; (ii) $B\cap \mathcal{Y}_{s}$ is 
of pure codimension 1; (iii) $\dim (B\cap \mathcal{Y}_{s})=n$. 
Condition (iii) holds iff the 
fiber of $q|_{B}:B\to S$ over $s$ has dimension $\geq n$. Since 
$q|_{B}$ is proper these points form  a closed subset of $S$ 
(\cite{Pe}~Corollary~IV.3.9). Its complement is the set satisfying (e).
\par
(f) Possibilities (ii) or (iii) occur iff the fiber 
of $q|_{B}:B\to S$ over $s$ has dimension $\geq n-1$. These points 
form a closed subset in $S$ whose complement is the set of points
satisfying (f).
\par
(g) and (h) Here we consider the restriction of $q$ to the closed subscheme 
$\mathcal{B}$. By hypothesis $q|_{\mathcal{B}}: \mathcal{B} \lto S$ is surjective.
By the arguments in the proof of (e), the condition that the fiber $\mathcal{B}_s$
is of pure codimension one in $\mathcal{Y}_{s}$ is equivalent to Condition~(e). Replacing 
$S$ by the corresponding open subscheme we may thus assume that every fiber of 
$q|_{\mathcal{B}}: \mathcal{B} \lto S$ is equidimensional of dimension $n-1$.
By \cite{Gr} Th\'{e}or\`{e}me~12.1.1 and the properness of $q:\mathcal{Y}\lto S$ the set of
points $s\in S$ such that $\mathcal{Y}_{s}$ has no embedded components is open in $S$. Replacing
$S$ by the corresponding open subscheme we may thus assume that this condition holds for every
fiber $\mathcal{Y}_{s}$. We claim that the scheme $\mathcal{B}$
is flat over $S$. Let $y\in \mathcal{B}$ and let $s=q(y)$.
Let $\mathcal{O}_y=\mathcal{O}_{\mathcal{Y},y}$ and let $\mathcal{O}_s=\mathcal{O}_{\mathcal{S},s}$.  Let $d_y\in \mathcal{O}_y$ 
be the germ of the discriminant. The restriction $d_y\otimes k(s)\in \mathcal{O}_y/\mathfrak{m}_s\mathcal{O}_y$ is not
a zero divisor, since the associated primes of $\mathcal{O}_y/\mathfrak{m}_s\mathcal{O}_y \cong \mathcal{O}_{_{\mathcal{Y}_{s}},y}$
are all of height 0, while every minimal prime ideal of $\mathcal{O}_{_{\mathcal{Y}_{s}},y}$ which contains $d_y\otimes k(s)$ is
of height 1, since $\mathcal{B}_s$ is of pure codimension 1 in $\mathcal{Y}_{s}$.
By \cite{Ma}~(20.E) we conclude that $\mathcal{O}_s \lto \mathcal{O}_y/d_y\mathcal{O}_y = \mathcal{O}_{\mathcal{B},y}$ is flat.
This proves that $q|_{\mathcal{B}}:\mathcal{B}\to S$ is flat. 
We may now apply to $q|_{\mathcal{B}}$ the arguments of (a) and (b).
\end{proof}
\begin{block}\label{s5.70}
Let $q:\mathcal{Y} \to \mathcal{Z}$ be a proper, flat morphism of 
schemes. Suppose $\mathcal{Z}$ is reduced and connected. Let 
$\mathcal{G}$ be a coherent locally free sheaf on $\mathcal{Y}$. Suppose 
$h^{0}(\mathcal{Y}_{s},\mathcal{G}_{s})$ is independent of $z\in 
\mathcal{Z}$ and nonzero. 
By Grauert's theorem (\cite{Mu} Ch.II,\S~5,Corollary~2) 
$\mathcal{H}=q_{*}\mathcal{G}$ is a coherent locally free sheaf. Let 
us denote by $\mathbb{G}\to \mathcal{Y}$, \; $f:\mathbb{H}\to 
\mathcal{Z}$ the vector bundles associated with $\mathcal{G}$ and 
$\mathcal{H}$ respectively: 
$\mathcal{G}=\mathcal{O}_{\mathcal{Y}}(\mathbb{G})$, 
$\mathcal{H}=\mathcal{O}_{\mathcal{Z}}(\mathbb{H})$, 
$\mathbb{H}_{z}=H^0(\mathcal{Y}_{z},\mathcal{G}_{z})$. The 
tautological morphism $q^{*}q_{*}\mathcal{G}\to \mathcal{G}$ 
corresponds to a vector bundle morphism 
\begin{equation}\label{es5.70}
\xymatrix{
\mathcal{Y}\times_{\mathcal{Z}}\mathbb{H}\ar[rr]\ar[rd]_{p_1}&&
\mathbb{G}\ar[dl]\\
&\mathcal{Y}
}
\end{equation}
It induces 
a section 
\begin{equation*}
\Psi :\mathcal{Y}\times_{\mathcal{Z}}\mathbb{H}
\lto \mathbb{G}\times_{\mathcal{Y}}(\mathcal{Y}\times_{\mathcal{Z}}
\mathbb{H}) = p_{1}^{*}\mathbb{G}
\end{equation*}
of the pull-back of $\mathbb{G}$.  If $(y,\eta)\in 
\mathcal{Y}\times_{\mathcal{Z}}\mathbb{H}$, then
$\Psi(y,\eta)=(\eta(y),(y,\eta))$.
\end{block}
\renewcommand{\theenumi}{\alph{enumi}}
\begin{lem}\label{s5.71}
Let $q:\mathcal{Y}\to \mathcal{Z}$ be a flat, proper morphism, where 
$\mathcal{Z}$ is a reduced, connected scheme. Suppose there is an integer 
$n\geq 1$ such that every irreducible component of every fiber of $q$
has dimension $n$. Suppose we have one of the following data:
\par 
($d=3$) Locally free sheaves $\mathcal{E}$ and $\mathcal{L}$ of 
$\mathcal{Y}$ of ranks 2 and 1 respectively. Let 
$\mathcal{G}=S^{3}\mathcal{E}\otimes \mathcal{L}$;
\par
($d=4$) 
Locally free sheaves $\mathcal{F}$, $\mathcal{E}$ and $\mathcal{L}$ of 
$\mathcal{Y}$ of ranks 2, 3 and 1 respectively. Let 
$\mathcal{G}=\Check{\mathcal{F}}\otimes S^{2}\mathcal{E}\otimes 
\mathcal{L}$;
\par
($d=5$) 
Locally free sheaves $\mathcal{F}$, $\mathcal{E}$ and $\mathcal{L}$ of 
$\mathcal{Y}$ of ranks 5, 4 and 1 respectively. Let 
$\mathcal{G}=\bigwedge^{2}{\mathcal{F}}\otimes \mathcal{E}\otimes 
\mathcal{L}$.
\par
Suppose $h^0(\mathcal{Y}_{z},\mathcal{G}_{z})$ is independent of $z$ 
and is nonzero. Consider the locally free sheaf 
$\mathcal{H}=q_{*}\mathcal{G}$ and let $f:\mathbb{H}\to \mathcal{Z}$ 
be the associated vector bundle with fibers $\mathbb{H}_{z}=
H^0(\mathcal{Y}_{z},\mathcal{G}_{z})$. Each of the following 
conditions on $\eta \in \mathbb{H}$ determines and open (possibly 
empty) subset of $\mathbb{H}$:
\begin{enumerate}
\item
If $f(\eta)=z$, then $\eta$ is of the right codimension for every 
$y\in \mathcal{Y}_{z}$;
\item
Assuming (a), if $\pi_{\eta}:X_{\eta}\to \mathcal{Y}_{z}$, $X_{\eta}
\subset \mathbf{P}(\mathcal{E}_{z})$, is the Gorenstein covering of 
degree $d$ associated with $\eta$, then $X_{\eta}$ is smooth;
\item
Assuming (a), the cover $X_{\eta}$ is irreducible and generically 
reduced;
\item
Assuming (a), the covering $\pi_{\eta}:X_{\eta}\to \mathcal{Y}_{z}$
is generically unramified (i.e. unramified over a dense open subset of 
$\mathcal{Y}_{z}$). 
\item
Assuming (a), the covering $\pi_{\eta}:X_{\eta}\to 
\mathcal{Y}_{z}$ is \'{e}tale;
\item
Assuming (a), in case (e) determines an empty set, 
$\mathcal{Y}_{z}$ has no embedded components and the 
discriminant subscheme of $\pi_{\eta}:X_{\eta}\to \mathcal{Y}_{z}$ is 
smooth of pure codimension one;
\item
Assuming (a), in case (e) determines an empty set,
$\mathcal{Y}_{z}$ has no embedded components and the 
discriminant subscheme of $\pi_{\eta}:X_{\eta}\to \mathcal{Y}_{z}$ is 
of codimension one, irreducible and generically reduced.
\end{enumerate}
Let $\mathbb{H}'\subset \mathbb{H}$ be the open subset defined by (a). 
Suppose $\mathbb{H}'\neq \emptyset$. Let $\mathcal{Y}'=\mathcal{Y}
\times_{\mathcal{Z}}\mathbb{H}'$. There exists a Gorenstein covering 
$\pi':\mathcal{X}'\to \mathcal{Y}'$ of degree $d$ such that for every 
$\eta   \in   \mathbb{H}'$   the   covering    $\mathcal{X}'_{\eta}\to 
\mathcal{Y}'_{\eta}$      is      isomorphic      to      $X_{\eta}\to 
\mathcal{Y}_{f(\eta)}$     via     the      canonical      isomorphism 
$\mathcal{Y}'_{\eta}\cong      \mathcal{Y}_{f(\eta)}$. Performing a 
base change of
$\mathbb{H}'$ by the open subsets determined by (b),  (c),  (d),  (e), 
(f) and (g) one obtains respectively  flat families of coverings  which  
parameterize  
all coverings, possessing the specified properties,
 present in  the  family  $\mathcal{X}'\to 
\mathcal{Y}'\to \mathbb{H}'$. 
\end{lem}
\begin{proof}
We focus on the case $d=5$ and leave to the reader the cases 
$d=3$ and $d=4$, which were also treated in \cite{K1} Lemma~2.8 and 
\cite{K2} Lemma~2.14.
\par
Let $\Psi:\mathcal{Y}\times_{\mathcal{Z}}\mathbb{H}\to 
p_{1}^{*}\mathbb{G}$ be the vector bundle section constructed in 
\S~\ref{s5.70}. One has $\mathcal{G}=\bigwedge^{2}\mathcal{F}\otimes 
\mathcal{E}\otimes \mathcal{L}$, so $p_{1}^{*}\mathcal{G}=
\bigwedge^{2}\tilde{\mathcal{F}}\otimes \tilde{\mathcal{E}}\otimes 
\tilde{\mathcal{L}}$, where 
$\tilde{\mathcal{F}}=p_{1}^{*}\mathcal{F}$,
$\tilde{\mathcal{E}}=p_{1}^{*}\mathcal{E}$,
$\tilde{\mathcal{L}}=p_{1}^{*}\mathcal{L}$.
Consider $\rho: \mathbf{P}(\tilde{\mathcal{E}})\to \mathcal{Y}\times
_{\mathcal{Z}}\mathbb{H}$. The section corresponds to an element 
$\Phi \in H^0(\mathbf{P}(\tilde{\mathcal{E}}),
(\bigwedge^{2}\rho^{*}\tilde{\mathcal{F}})\otimes 
\rho^{*}\tilde{\mathcal{L}}(1))$. The set 
\begin{equation*}
\Gamma = \{x\in \mathbf{P}(\tilde{\mathcal{E}})|
\, \text{rank}\, \Phi(x)\leq 2\}
\end{equation*}
is closed in $\mathbf{P}(\tilde{\mathcal{E}})$, it is of codimension 
$\leq 3$ and its intersection with any fiber 
$\mathbf{P}(\tilde{\mathcal{E}})_{(y,\eta)}$ is nonempty. The 
projection $\Gamma \to \mathcal{Y}\times_{\mathcal{Z}}\mathbb{H}$
is therefore surjective and for every $(y,z)\in 
\mathcal{Y}\times_{\mathcal{Z}}\mathbb{H}$ there are two possibilities:
either $\dim \Gamma_{(y,z)}=0$, or $\dim \Gamma_{(y,z)}\geq 1$. Let us 
denote by $\Sigma$ the subset of 
$\mathcal{Y}\times_{\mathcal{Z}}\mathbb{H}$ where the second 
alternative holds. By the properness of $\Gamma \to 
\mathcal{Y}\times_{\mathcal{Z}}\mathbb{H}$ the set $\Sigma$ is closed 
(\cite{Pe}~Corollary~IV.3.9). By assumption $q:\mathcal{Y}\to \mathcal{Z}$ is 
proper, so $p_{2}: \mathcal{Y}\times_{\mathcal{Z}}\mathbb{H} \to 
\mathbb{H}$ is proper as well and therefore $p_{2}(\Sigma)$ is closed 
in $\mathbb{H}$. The subset $\mathbb{H}'=\mathbb{H}\setminus 
p_{2}(\Sigma)$, consisting of points where (a) holds, is 
therefore open in $\mathbb{H}$.
\par
Let $\mathcal{F}'$, $\mathcal{E}'$, $\mathcal{L}'$ be the restrictions 
of $\tilde{\mathcal{F}}$,
$\tilde{\mathcal{E}}$, $\tilde{\mathcal{L}}$ on 
$\mathcal{Y}'=\mathcal{Y}\times_{\mathcal{Z}}\mathbb{H}'$. Restricting $\Psi$ 
one
obtains a section $\Psi'\in H^0(\mathcal{Y}',
\bigwedge^2\mathcal{F}'\otimes \mathcal{E}'\otimes \mathcal{L}')$, which 
has the right codimension at every $(y,\eta)\in \mathcal{Y}'$. Let
$\mathcal{X}'\subset \mathbf{P}(\mathcal{E}')$ be the 
associated
Pfaffian 
subscheme (cf. \S~\ref{s5.80}). The projection $\mathcal{X}'\to 
\mathcal{Y}'$ is a Gorenstein covering of degree 5 according to 
Proposition~\ref{s5.81}. Consider the resolution of 
$\mathcal{O}_{\mathcal{X}'}$ given in \eqref{es5.81}. Let $\eta \in
\mathbb{H}'$ and let $f(\eta)=z$. The base change $\mathcal{Y}_{z}
\times \{\eta\}\to \mathcal{Y}\times_{\mathcal{Z}}\mathbb{H}'$ 
transforms the above resolution into a resolution of 
$\mathcal{O}_{\mathcal{X}'_{\eta}}$, since 
$\mathcal{O}_{\mathcal{X}'}$ is $\mathcal{O}_{\mathcal{Y}'}$-flat
(see Lemma~\ref{s1.p2}). The 
functoriality of the Buchsbaum-Eisenbud complex \eqref{es5.82} 
implies that $\mathcal{X}'_{\eta}\cong X_{\eta}\subset 
\mathbf{P}(\mathcal{E}\otimes \mathcal{O}_{\mathcal{Y}_{z}})$.
\par
Now, in order to prove the openness of conditions (b), (c), (d), 
(e), (f) and (g) one applies Proposition~\ref{s5.69} to the family of 
coverings $\mathcal{X}'\to \mathcal{Y}'\to \mathbb{H}'$.
\par
The last statement of the lemma is clear.
\end{proof}
Let $A$ be a Noetherian ring, $S=A[x_{0},x_{1},\ldots,x_{r}]$, $P = 
Proj\, S=\mathbb{P}^{r}_{A}$. Let $\mathcal{G}$ be a quasicoherent 
sheaf on $P$. We denote by $H^{i}_{*}(\mathcal{G})$ the graded 
$S$-module $H^{i}_{*}(\mathcal{G}):=
\oplus_{n\in \mathbb{Z}}H^{i}(P,\mathcal{G}(n))$.  Recall that 
$\mathcal{G}$ is the sheafification of 
$G=H^{0}_{*}(\mathcal{G})$:\; $\mathcal{G}=G^{\sim}$ (cf. \cite{Ha} Ch.II, 
\S~5).
\renewcommand{\theenumi}{\roman{enumi}}
\begin{lem}\label{s5.dop}
Suppose $\mathcal{G}$ is a coherent sheaf on $P=\mathbb{P}^{r}_{A}$ which has 
a 
resolution of the type:
\begin{equation*}
\mathcal{F}_{\bullet}: 0\to \mathcal{F}_{m}\to \ldots \to 
\mathcal{F}_{i}\to \ldots \to \mathcal{F}_{0}\to \mathcal{G}\to 0
\end{equation*}
where $m\leq r$ and $\mathcal{F}_{i}\cong \oplus _{j=1}^{\beta_{i}}
\mathcal{O}_{P}(-a_{ij})$, $a_{ij}\in \mathbb{Z}$, $\beta_{j}\in 
\mathbb{N}$ for $\forall i,j$. Then:
\begin{enumerate}
\item
$H^{i}_{*}(\mathcal{G})=0$ for $0<i<r-m$;
\item
$G=H^{0}_{*}(\mathcal{G})$ is a finitely generated graded $S$-module;
\item
$\mathcal{F}_{\bullet}$ is the sheafification of the exact complex of 
graded $S$-modules
\begin{equation*}
H^{0}_{*}(\mathcal{F}_{\bullet}): 0\to F_{m}\to \ldots \to F_{i}\to
\ldots \to F_{0}\to G\to 0
\end{equation*}
where $F_{i}=\oplus_{j=1}^{\beta_{i}}S(-a_{ij})$.
\end{enumerate}
\end{lem}
\begin{proof}
The lemma is proved by induction on $m$. If $m=0$ the statements 
follow from Serre's theorem (see \cite{Ha} Ch. III, \S~5). Suppose 
$m\geq 1$. Let $\mathcal{K}_{0}=Ker(\mathcal{F}_{0}\to \mathcal{G})$. 
By the induction assumption the lemma holds for $\mathcal{K}_{0}$ with 
$m$ replaced by $m-1$. The long exact sequence of cohomology 
associated with $0\to \mathcal{K}_{0}\to \mathcal{F}_{0}\to 
\mathcal{G}\to 0$ yields $H^{i}_{*}(\mathcal{G})\cong 
H^{i+1}_{*}(\mathcal{K}_{0})=0$ for $0<i<r-m$ and the 
exactness of $0\to H^{0}_{*}(\mathcal{K}_{0})\to 
H^{0}_{*}(\mathcal{F}_{0})\to H^{0}_{*}(\mathcal{G})\to 0$, since 
$H^{1}_{*}(\mathcal{K}_{0})=0$. This implies (ii) and (iii).
\end{proof}
Our next goal is to study the action of 
$Aut(F)\times Aut(E)$ on $H^0(Y,\bigwedge^2F\otimes E\otimes 
(\det E)^{-1})$ in connection with the relation of
equivalence of coverings of $Y$.
\begin{lem}\label{s5.84}
Let $Y$ be an irreducible variety.
Let $E$, $F$ be locally free 
sheaves on $Y$ of ranks 4, 5 respectively such that $\det F\cong (\det E)^2$. 
Suppose 
$\eta, \eta' \in H^0(Y,\bigwedge^2F\otimes E\otimes (\det E)^{-1})$ have 
the right codimension at every $y\in Y$. Suppose $X_{\eta} = X_{\eta'}
\subset \mathbf{P}(E)$. Then there exists an element $\sigma \in 
Aut(F)$ such that $(\wedge^{2}\sigma \otimes id)\eta = h \, \eta'$ 
for 
some invertible element $h \in \Gamma(Y,\mathcal{O}_{Y})^{*}$. 
Furthermore $\sigma$ is unique up to multiplication by an element of
$\Gamma(Y,\mathcal{O}_{Y})^{*}$.
\end{lem}
\begin{proof}
Let $X=X_{\eta}=X_{\eta'}$. The resolution of $\mathcal{O}_{X}$ (see 
\eqref{es5.81a}) is unique up  to unique isomorphism according to 
\cite{CE}~Theorem~2.1. We thus obtain a commutative diagram
\begin{equation} \label{es5.84}
\xymatrix{
0\ar[r]&\rho^{*}(F_{3})(-5)\ar[d]^{\alpha_{3}}\ar[r]^{^{t}g}&
\rho^{*}(F_{2})(-3)\ar[d]^{\alpha_{2}}\ar[r]^{f}&
\rho^{*}(F_{1})(-2)\ar[d]^{\alpha_{1}}\ar[r]^-{g}&
\mathcal{O}_{\mathbf{P}(E)}\ar[d]^{\alpha_{0}=id}\\
0\ar[r]&\rho^{*}(F_{3})(-5)\ar[r]^{^{t}g'}&
\rho^{*}(F_{2})(-3)\ar[r]^{f'}&
\rho^{*}(F_{1})(-2)\ar[r]^-{g'}&
\mathcal{O}_{\mathbf{P}(E)}
}
\end{equation}
where $F_{1}=F$, $F_{2}=\Check{F}\otimes \det E$, $F_{3}=\det E$. 
Here the vertical morphisms are isomorphisms and $f, g,\; {^{t}g}$ 
(resp. $f', g',\; {^{t}g')}$ are obtained from $\eta$ (resp. $\eta'$) as 
in \eqref{es5.82a}. Since $rk(F_{3})=1$ one has that 
$\alpha_{3}=c\cdot id$, where $c\in \Gamma(Y,\mathcal{O}_{Y})^{*}$. 
The resolutions are self-dual. Applying the functor 
$\mathcal{H}om_{\mathbf{P}(E)}(-,\rho^{*}(F_{3})(-5))$ and using that 
${^{t}f}=-f$, ${^{t}f'}=-f'$ one obtains the commutative diagram
\begin{equation}\label{es5.85}
\xymatrix{
0\ar[r]&\rho^{*}(F_{3})(-5)\ar[d]^{^{t}\alpha_{0}^{-1}}\ar[r]^{^{t}g}&
\rho^{*}(F_{2})(-3)\ar[d]^{^{t}\alpha_{1}^{-1}}\ar[r]^{f}&
\rho^{*}(F_{1})(-2)\ar[d]^{^{t}\alpha_{2}^{-1}}\ar[r]^-{g}&
\mathcal{O}_{\mathbf{P}(E)}\ar[d]^{c^{-1}\cdot id}\\
0\ar[r]&\rho^{*}(F_{3})(-5)\ar[r]^{^{t}g'}&
\rho^{*}(F_{2})(-3)\ar[r]^{f'}&
\rho^{*}(F_{1})(-2)\ar[r]^-{g'}&
\mathcal{O}_{\mathbf{P}(E)}
}
\end{equation}
The chain maps $\{c^{-1}\alpha_{i}\}$ and $\{^{t}\alpha_{3-i}^{-1}\}$ 
are liftings of the morphism $c^{-1}\cdot 
id:\mathcal{O}_{\mathbf{P}(E)}\to \mathcal{O}_{\mathbf{P}(E)}$. We 
claim $^{t}\alpha_{3-i}^{-1} = c^{-1}\alpha_{i}$ for every 
$i$. Let us cover $Y$ by affine open sets 
$Y=\cup_{\lambda}U_{\lambda}$ such that $E|_{U_{\lambda}}$ and 
$F|_{U_{\lambda}}$ are trivial for every $\lambda$. It suffices to 
verify the equalities $^{t}\alpha_{3-i}^{-1} = c^{-1}\alpha_{i}$
restricting to every $\rho^{-1}U_{\lambda}$. The restrictions of 
\eqref{es5.84} and \eqref{es5.85} to these open subsets are 
sheafifications of minimal resolutions of  graded modules
by Lemma~\ref{s5.dop}
and we obtain the claim, since the possible homotopies are zero. Now 
we have 
\begin{equation}\label{es5.85a}
f' = \alpha_{1}\circ f \circ \alpha_{2}^{-1} = c\; 
{^{t}\alpha}_{2}^{-1}\circ f \circ \alpha_{2}^{-1}.
\end{equation}
The morphism $\alpha_{2}^{-1}$ is induced by an element of $Aut(F)$. 
Looking at the definition of the differentials $f$ and $f'$ in the 
Buchsbaum-Eisenbud complex we see that \eqref{es5.85a} is equivalent 
to the equality $\wedge^{2}\sigma \otimes id\, (\eta) = h \, \eta'$ 
for some $\sigma \in Aut(F)$ and $h \in 
\Gamma(Y,\mathcal{O}_{Y})^{*}$.
\par
In order to prove the uniqueness statement it suffices to show that for 
$\sigma \in Aut(F)$ and $\eta \in H^0(Y,\bigwedge^2F\otimes E\otimes 
(\det E)^{-1})$ of right codimension at every $y\in Y$ the equality 
$(\wedge^{2}\sigma \otimes id)\eta = h \, \eta$, $h \in 
\Gamma(Y,\mathcal{O}_{Y})^{*}$ implies $\sigma = c\cdot id_{F}$, $c\in 
\Gamma(Y,\mathcal{O}_{Y})^{*}$. Let us first consider the case when 
$Y$ is one point. Then $\eta \in \bigwedge^{2}k^{5}\otimes k^{4}$ and 
$\sigma \in Aut(k^{5})$. Multiplying $\sigma$ by a constant we may 
assume $(\wedge^{2}\sigma \otimes id)\eta = \eta$. Here $X\subset 
\mathbb{P}^{3}$ is a nondegenerate, arithmetically Gorenstein 
subscheme with $\dim X=0$, $\deg X=5$. Let $S=k[x_{0},\ldots,x_{3}]$. 
The linear map
$\sigma$ induces the following 
automorphism of the graded Buchsbaum-Eisenbud complex associated with $\eta$. 
\begin{equation*}
\xymatrix{
0\ar[r]&S(-5)\ar[d]^{c^{-1}id}\ar[r]^{^{t}g}&
S(-3)\ar[d]_{\sigma}\ar[r]^{f}&
S(-2)\ar[d]^{^{t}\sigma^{-1}}\ar[r]^-{g}&
S\ar[d]^{c\cdot id}\ar[r]&S_{X}\ar[r]\ar[d]^{c\cdot id}&0\\
0\ar[r]&S(-5)\ar[r]^{^{t}g}&
S(-3)\ar[r]^{f}&
S(-2)\ar[r]^-{g}&
S\ar[r]&S_{X}\ar[r]&0
}
\end{equation*}
Here $c=\det(\sigma)^{-1}$. The homomorphism $c\cdot id:S\to S$ has a 
unique lifting to a chain map of the acyclic complex since the 
possible homotopies are 0. Therefore $\sigma = c\cdot id$, 
$c^{-1}\cdot id=c\cdot id$. This shows that $c=\pm 1$, $\sigma = \pm 
id$. We conclude that, when $Y$ is one point, if $(\wedge^{2}\sigma
\otimes id)\eta = a\eta$, $a\in k^{*}$, then $\sigma = \pm 
\sqrt{a}\cdot id$. 
\par
Suppose now $Y$ is an arbitrary irreducible variety. Let 
$(\wedge^{2}\sigma \otimes id)\eta = h \eta$ for some $h 
\in \Gamma(Y,\mathcal{O}_{Y})^{*}$. Evaluating this equality at every 
$y\in Y$ we obtain that $\sigma(y)\in k^{*}\cdot id$ for $\forall y\in 
Y$. One has 
\begin{equation*}
Aut(F) \subset H^0(Y,End(Y)) = H^0(Y,\mathcal{O}_{Y})\cdot id_{F} \oplus 
H^0(Y,ad'(F))
\end{equation*}
Therefore $\sigma \in \Gamma(Y,\mathcal{O}_{Y})^{*}\cdot id_{F}$.
\end{proof}
In the next lemmas we use the notation 
$k[Y]:=\Gamma(Y,\mathcal{O}_{Y})$.
\begin{lem}[Case $d=5$]\label{s5.87a}
Let $Y$ be an irreducible variety.
Let $E$, $F$ be locally free 
sheaves on $Y$ of ranks 4 and 5 respectively 
such that $\det F\cong (\det E)^2$.
\par
(i) Let $\eta,\eta' \in H^0(Y,\bigwedge^2F\otimes E\otimes (\det E)^{-1})$
be of the right codimension at every $y\in Y$. Then the covering 
$\pi_{\eta}:X_{\eta}\to Y$ is equivalent to $\pi_{\eta'}:X_{\eta'}\to Y$
if and only if $(\wedge^{2}\sigma \otimes \varphi \otimes id)\eta = 
h \eta'$ for some $\sigma \in Aut(F)$, $\varphi \in Aut(E)$, $h\in 
k[Y]^{*}$.
\par
(ii) Suppose furthermore $\pi_{\eta}:X_{\eta}\to Y$ is generically 
unramified. The pairs $(\sigma,\varphi)$ such that \linebreak
$(\wedge^{2}\sigma 
\otimes \varphi \otimes id)\eta = h \eta$ for some $h\in k[Y]^{*}$ 
form a subgroup in $Aut(F)\times Aut(E)$ whose quotient by 
the normal subgroup
$k[Y]^{*}id_{F}\times k[Y]^{*}id_{E}$ is a finite group isomorphic to 
$Aut(X_{\eta}/Y)$.
\par
(iii) Suppose moreover that 
$\pi_{\eta}:X_{\eta}\to Y$ is generically 
unramified and $Aut(X_{\eta}/Y)=\{1\}$. Then $(\sigma,\varphi)$ 
satisfies $(\wedge^{2}\sigma \otimes \varphi \otimes id)\eta = h 
\eta$ for some $h\in k[Y]^{*}$ if and only if $\sigma =a\cdot id_{F}$, 
$\varphi =b\cdot id_{E}$ for some $a,b\in k[Y]^{*}$.
\end{lem}
\begin{proof}
(i) Let $\pi = \pi_{\eta}$, $\pi'=\pi_{\eta'}$. If $(\wedge^{2}\sigma \otimes 
\varphi \otimes id)\eta = h 
\eta'$, then $X_{\eta'}\cong {^{t}\varphi}^{-1}(X_{\eta})$, so 
$X_{\eta}\to Y$ is equivalent to $X_{\eta'}\to Y$. Vice versa, suppose 
there is an isomorphism $f$ such that the following diagram commutes
\begin{equation*}
\xymatrix{
X_{\eta}\ar[rr]^{f}\ar[rd]_{\pi}&&
X_{\eta'}\ar[dl]^{\pi'}\\
&Y
}
\end{equation*}
From the functoriality of the relative dualizing sheaf and the trace 
map one has a commutative diagram
\begin{equation*}
\xymatrix{
0\ar[r]&E\ar[d]^{\varphi}\ar[r]&
\pi'_{*}(\omega_{X_{\eta'}/Y})\ar[d]^{\cong}\ar[r]^-{t_{\pi'}}&
\mathcal{O}_{Y}\ar[d]^{=}\ar[r]&
0\\
0\ar[r]&E\ar[r]&
\pi_{*}(\omega_{X_{\eta}/Y})\ar[r]^-{t_{\pi}}&
\mathcal{O}_{Y}\ar[r]&
0
}
\end{equation*}
It induces the commutative diagram
\begin{equation*}
\xymatrix{
X_{\eta}\ar@{^{(}->}[r]\ar[d]_{f}&\mathbf{P}(E)\ar[d]^{{^{t}\varphi}
^{-1}}\\
X_{\eta'}\ar@{^{(}->}[r]&\mathbf{P}(E)
}
\end{equation*}
where the  embeddings are induced by the compositions
\begin{equation*}
\pi^{*}E \to \pi^{*}\pi_{*}(\omega_{X_{\eta}/Y})\to \omega_{X_{\eta}/Y}, 
\quad
\pi'^{*}E \to \pi'^{*}\pi'_{*}(\omega_{X_{\eta'}/Y})\to \omega_{X_{\eta'}/Y}
\end{equation*}
(cf. \cite{CE} Theorem~2.1~(ii)). This shows that 
\begin{equation*}
X_{\eta'} = {^{t}\varphi}^{-1}(X_{\eta}) = X_{id\otimes \varphi 
\otimes id (\eta)}.
\end{equation*}
According to Lemma~\ref{s5.84} there exists $\sigma \in Aut(F)$ such 
that $(\wedge^{2}\sigma \otimes \varphi \otimes id)\eta = h\eta'$ for 
some $h\in k[Y]^{*}$.
\par
(ii) $Aut(E)/k[Y]^{*}id_E$ acts on $\mathbf{P}(E)$ by the formula 
$\overline{\varphi}\cdot x = {^{t}\varphi}^{-1}(x)$. Let 
$G\subset Aut(E)/k[Y]^{*}id_E$ be the subgroup which leaves invariant 
$X_{\eta}$. The assumption that $X_{\eta}\to Y$ is generically 
unramified and Part~(i) imply that $G\cong Aut(X_{\eta}/Y)$. The pairs 
$(\sigma,\varphi)$ such that $(\wedge^{2}\sigma \otimes \varphi 
\otimes id)\eta = h\eta$ for some $h\in k[Y]^{*}$ form a subgroup of 
$Aut(F)\times Aut(E)$. Its quotient by the normal subgroup 
$k[Y]^{*}\cdot id_{F}\times k[Y]^{*}\cdot id_{E}$ projects bijectively 
to $G$ according to Lemma~\ref{s5.84}. This proves (ii).
\par
(iii) This is immediate from (ii) since $G\cong Aut(X_{\eta}/Y)\cong 
\{1\}$. 
\end{proof}
With similar and simpler arguments, left to the reader, one can prove 
analogs of Lemma~\ref{s5.87a} in the cases of coverings of degree $d=3$ and 
$d=4$. 
\begin{lem}[Case $d=3$]\label{s5.89}
Let $Y$ be an irreducible variety.
Let $E$ be a locally free 
sheaf on $Y$ of rank 2.
\par
(i) Let $\eta,\eta' \in H^0(Y,S^3E\otimes (\det E)^{-1})$ be of the 
right codimension at every $y\in Y$. Then the covering 
$\pi_{\eta}:X_{\eta}\to Y$ is equivalent to $\pi_{\eta'}:X_{\eta'}\to 
Y$ if and only if $(S^{3}\varphi \otimes id)\eta = h\eta'$ for some 
$\varphi \in Aut(E)$ and $h\in k[Y]^{*}$. 
\par
(ii) Suppose furthermore $\pi_{\eta}:X_{\eta}\to Y$ is generically
unramified. The automorphisms $\varphi \in Aut(E)$ such that 
$(S^{3}\varphi \otimes id)\eta = h\eta$ for some $h\in k[Y]^{*}$ form 
a subgroup in $Aut(E)$ whose quotient by the normal subgroup 
$k[Y]^{*}id_{E}$ is isomorphic to $Aut(X_{\eta}/Y)$.
\par
(iii) Suppose moreover 
$\pi_{\eta}:X_{\eta}\to Y$ is generically
unramified and $Aut(X_{\eta}/Y)=\{1\}$. Then 
$(S^{3}\varphi \otimes id)\eta = h\eta$ for some $h\in k[Y]^{*}$ if 
and only if $\varphi = a\cdot id_{E}$ for some $a\in k[Y]^{*}$. 
\end{lem}
\begin{lem}[Case $d=4$]\label{s5.89a}
Let $Y$ be an irreducible variety.
Let $E$, $F$ be locally free 
sheaves on $Y$ of ranks 3 and 2 respectively such that $\det F\cong \det E$.
\par
(i) Let $\eta,\eta' \in H^0(Y,\Check{F}\otimes S^2E)$ be of the right 
codimension at every $y\in Y$. 
Then the covering 
$\pi_{\eta}:X_{\eta}\to Y$ is equivalent to $\pi_{\eta'}:X_{\eta'}\to 
Y$ if and only if $({^t\sigma}\otimes S^{2}\varphi)\eta = h\eta'$
for some $\sigma \in Aut(F)$, $\varphi \in Aut(E)$ and $h\in k[Y]^{*}$. 
\par
(ii) Suppose furthermore $\pi_{\eta}:X_{\eta}\to Y$ is generically 
unramified. The pairs $(\sigma,\eta)$ such that 
$({^t\sigma}\otimes S^{2}\varphi)\eta = h\eta$ for some 
$h\in k[Y]^{*}$ form a subgroup of $Aut(F)\times Aut(E)$ whose 
quotient by $k[Y]^{*}id_{F}\times k[Y]^{*}id_{E}$ is a finite group 
isomorphic to $Aut(X_{\eta}/Y)$. 
\par
(iii) Suppose moreover $Aut(X_{\eta}/Y)=\{1\}$. Then $(\sigma,\varphi)$ 
satisfies $({^t\sigma}\times S^{2}\varphi)\eta = h\eta$ for some 
$h\in k[Y]^{*}$ if and only if $\sigma = a\cdot id_{F}$, $\varphi = 
b\cdot id_{E}$ for some $a,b\in k[Y]^{*}$.
\end{lem}
\begin{pro}\label{s5.89b}
Let $Y$ be a complete irreducible variety. Let $d=3,\, 4$ or $5$.
If 
$d=3$ let $E$ be a locally free sheaf on $Y$ of rank 2. Let $R_{3}=
S^3E\otimes (\det E)^{-1}$. If $d=4$ let $E$ and $F$ be locally free 
sheaves on $Y$ of ranks 3 and 2 respectively, such that $\det E\cong 
\det F$. Let $R_{4}=\Check{F}\otimes S^2E$. If $d=5$ let $E$ and $F$ 
be locally free sheaves on $Y$ of ranks 4 and 5 respectively, such 
that $(\det E)^{\otimes 2}\cong \det F$. Let $R_{5}=\bigwedge^2F\otimes 
E\otimes (\det E)^{-1}$. The subset $W\subset \mathbb{P}H^0(Y,R_{d})$,
consisting of $\langle \eta\rangle$ such that $\eta$ is of the right 
codimension at every $y\in Y$ and $\pi_{\eta}:X_{\eta}\to Y$ is 
generically unramified, is open. Suppose $W\neq \emptyset$. Then 
the group $Aut(E)/k^{*}$ (when $d=3$), resp. $Aut(F)/k^{*}\times Aut(E)/k^{*}$ 
(when $d=4$ or $5$) acts with finite stabilizers on $W$. The action is 
free on the subset of $W$ consisting of  $\langle \eta \rangle$ 
which satisfy $Aut(X_{\eta}/Y)=\{1\}$. Two elements $\langle \eta 
\rangle, \langle \eta' \rangle \in W$  belong to the same orbit if and 
only if $\pi_{\eta}:X_{\eta}\to Y$ is equivalent to 
$\pi_{\eta'}:X_{\eta'}\to Y$.
\end{pro}
\begin{proof}
The set $W\subset \mathbb{P}H^0(Y,R_{d})$ is open according to 
Lemma~\ref{s5.71} (here $q:\mathcal{Y}\to \mathcal{Z}$ is $Y\to 
\Spec k$). The other statements follow from Lemma~\ref{s5.89}. 
Lemma~\ref{s5.89a} and Lemma~\ref{s5.87a}.
\end{proof}
We refer to \cite{CE} Theorem~3.6 and Theorem~4.5, and to 
\cite{C2} Theorem~4.4 for the following result of Bertini type
(see also Lemma~\ref{s5.71}).
\begin{pro}\label{s1.7}
Let the hypothesis and notation be as in Proposition~\ref{s5.89b}. Suppose 
furthermore that: 
$k=\mathbb{C}$; the irreducible variety  $Y$ is smooth and projective; $\dim 
Y\leq 3$ if $d=3$
or $4$ and $\dim Y\leq 2$ if $d=5$. 
 Suppose $R_{d}$ is generated by its global 
sections $H^0(Y,R_{d})$. Then the following two subsets of 
$H^0(Y,R_{d})$ are nonempty and open:
\begin{align*}
H_{rc} &=\{\eta \, | \, \eta \:\text{has the right codimension at 
every}\: y\in 
Y\}\\
H_{s}  &=\{\eta \, | \, \eta \in H_{rc},\; X_{\eta}\: \text{is smooth}\}
\end{align*}
Furthermore if $\eta\in H_{rc}$ (resp. if $\eta\in H_{s}$), then $X_{\eta}$
 is 
connected (resp. irreducible) if and 
only if $H^0(Y,\Check{E})=0$.
\end{pro}

\section{Construction of families of coverings with rational parameter 
varieties}\label{s2}
We fix a smooth, irreducible, projective curve $Y$ of genus $g\geq 1$. 
\renewcommand{\theenumi}{\roman{enumi}}
\begin{lem}\label{s2.20}
Let $e\in \mathbb{Z}$. There exists an irreducible, nonsingular 
variety $S$ and a locally free sheaf $\mathcal{E}$ of rank 2 
over $S\times Y$ with the following properties.
\begin{enumerate}
\item
$\mathcal{E}_{s}$ has degree $e$, it is stable if $g\geq 2$, and  regular 
polystable if $g=1$, for every $s\in S$.
\item
Every stable/regular polystable locally free sheaf of rank 2 and 
degree $e$ over $Y$ is isomorphic to $\mathcal{E}_{s}$ for some $s\in 
S$. 
\item
The morphism $\delta:S\to Pic^{e}Y$, given by 
$\delta(s)=\det(\mathcal{E}_{s})$, is surjective.
\item
The fiber $S_{A}=\delta^{-1}(A)$ is nonempty, irreducible, nonsingular 
and rational for every $A\in Pic^{e}Y$.
\item
For every $A\in Pic^{e}Y$, and every $s\in S_{A}$, the Kodaira-Spencer 
map $\kappa:T_{s}S\to H^1(Y,ad'(\mathcal{E}_{s}))$ is surjective, 
where $ad'(\mathcal{E}_{s})$ is the kernel of 
$Tr:End(\mathcal{E}_{s})\to \mathcal{O}_{Y}$.
\end{enumerate}
\end{lem}
\renewcommand{\theenumi}{\roman{enumi}}
\begin{proof}
We apply Proposition~\ref{s1.18a} with $T=\emptyset$ and 
we denote by $\tilde{S}$ the variety constructed in it.
Let $S\subset \tilde{S}$ be 
the Zariski open subset whose 
points correspond to the stable/regular polystable sheaves.
\end{proof}
\begin{lem}\label{s2.22}
Let $e\in \mathbb{Z}$. There exists an irreducible, nonsingular 
variety $S$ and locally free sheaves $\mathcal{E}$ and $\mathcal{F}$ 
over $S\times Y$ of ranks 3 and 2 respectively, such that $\det 
\mathcal{E}\cong \det \mathcal{F}$, and the following properties are 
satisfied.
\begin{enumerate}
\item
$\mathcal{E}_{s}$ and $\mathcal{F}_{s}$ have degree $e$, they are 
stable if $g\geq 2$, and  regular polystable if $g=1$,
for every $s\in S$. 
\item
Every pair $E,F$ of stable/regular polystable locally free sheaves 
over $Y$, such that $rk(E)=3,rk(F)=2, \deg E = \deg F = e$, and $\det 
E\cong \det F$, is isomorphic to $\mathcal{E}_{s},\mathcal{F}_{s}$ for 
some $s\in S$.
\item
The morphism $\delta:S\to Pic^{e}Y$, given by 
$\delta(s)=\det(\mathcal{E}_{s})$ is surjective.
\item
The fiber $S_{A}=\delta^{-1}(A)$ is nonempty, irreducible, 
nonsingular, and rational for every $A\in Pic^{e}Y$. 
\item 
For every $A\in Pic^{e}Y$ and every $s\in S_{A}$ the Kodaira-Spencer 
map 
\[
\kappa : T_{s}S_{A}\lto H^1(Y,ad'(\mathcal{E}_{s}))\oplus 
H^1(Y,ad'(\mathcal{F}_{s}))
\]
is surjective.
\end{enumerate}
\end{lem}
\begin{proof}
Let $y_{0}\in Y$. Let us choose $d=e+6N$ sufficiently large. Let $J$ 
be the Jacobian variety of $Y$. Let $\mathcal{P}_{d}$ be the 
Poincar\'{e} invertible sheaf on $J\times Y$, satisfying 
$\mathcal{P}_{d}|_{\{L\}\times Y}\cong L\otimes 
\mathcal{O}_{Y}(dy_{0})$ and 
$\mathcal{P}_{d}|_{J\times \{y_{0}\}}\cong
\mathcal{O}_{J}$. Using Proposition~\ref{s1.18a} one obtains a Zariski 
open subset $S_{1}\subset \mathbb{V}_1$, where 
$\pi_{1}:\mathbb{V}_{1}\to J$ is a vector bundle, and an extension
\[
0\to \mathcal{O}^{\oplus 2}_{S_{1}\times Y}\to \mathcal{E}'\to 
(\pi_{1}\times id)^{*}\mathcal{P}_{d} \to 0,
\]
such that $\mathcal{E}'_{s}$ is stable/regular polystable for every 
$s\in S_{1}$. This family includes all stable/regular polystable 
locally free sheaves of degree $d$ and rank 3. In a similar manner one 
treats the rank 2 case obtaining $\pi_{2}:\mathbb{V}_{2}\to J,\: 
S_{2}\subset \mathbb{V}_{2}$ and an extension
\[
0\to \mathcal{O}_{S_{2}\times Y}\to \mathcal{F}'\to 
(\pi_{2}\times id)^{*}\mathcal{P}_{d} \to 0.
\]
One has $\det \mathcal{E}'\cong (\pi_{1}\times 
id)^{*}\mathcal{P}_{d}$, $\det \mathcal{F}'\cong (\pi_{2}\times 
id)^{*}\mathcal{P}_{d}$. Let $S=S_{1}\times _{J}S_{2}\subset 
\mathbb{V}_{1}\oplus \mathbb{V}_{2}$ be the open inclusion. Consider 
the cartesian diagram 
\[
\xymatrix{
S\ar[r]^{\rho_{2}}\ar[d]_{\rho_{1}}&S_{2}\ar[d]^{\pi_{2}}\\
S_{1}\ar[r]^-{\pi_{1}}&J
}
\]
Let $\mathcal{E}=(\rho_{1}\times id)^{*}\mathcal{E}'\otimes 
p_{Y}^{*}\mathcal{O}_{Y}(-2Ny_{0})$ and let 
$\mathcal{F}=(\rho_{2}\times id)^{*}\mathcal{F}'\otimes 
p_{Y}^{*}\mathcal{O}_{Y}(-3Ny_{0})$. Then $S,\mathcal{E}$ and 
$\mathcal{F}$ satisfy the required properties of the lemma.
\end{proof}
\begin{lem}\label{s2.24}
Let $e\in \mathbb{Z}$. There exists an irreducible, nonsingular 
variety $S$ and locally free sheaves $\mathcal{E}$ and $\mathcal{F}$ 
over $S\times Y$ of ranks 4 and 5 respectively, such that $\det 
\mathcal{F}\cong (\det \mathcal{E})^{2}$, and the following properties are 
satisfied.
\begin{enumerate}
\item
$\mathcal{E}_{s}$ and $\mathcal{F}_{s}$ have degrees $e$ and $2e$ 
respectively, they are 
stable if $g\geq 2$, and  regular polystable if $g=1$,
for every $s\in S$. 
\item
Every pair $E,F$ of stable/regular polystable locally free sheaves 
over $Y$, such that $rk(E)=4,rk(F)=5, \deg E = e,\: \deg F = 2e$, 
and $\det F\cong (\det E)^{2}$, is isomorphic to 
$\mathcal{E}_{s},\mathcal{F}_{s}$ for 
some $s\in S$.
\item
The morphism $\delta:S\to Pic^{e}Y$, given by 
$\delta(s)=\det(\mathcal{E}_{s})$ is surjective.
\item
The fiber $S_{A}=\delta^{-1}(A)$ is nonempty, irreducible, 
nonsingular, and rational for every $A\in Pic^{e}Y$. 
\item 
For every $A\in Pic^{e}Y$ and every $s\in S_{A}$ the Kodaira-Spencer 
map 
\[
\kappa : T_{s}S_{A}\lto H^1(Y,ad'(\mathcal{E}_{s}))\oplus 
H^1(Y,ad'(\mathcal{F}_{s}))
\]
is surjective.
\end{enumerate}
\end{lem}
\begin{proof}
This is proved in the same way as Lemma~\ref{s2.22} with the following 
modifications. One lets $d=e+20N,\, N\gg 0$ and considers extensions
\begin{align*}
&0\to \mathcal{O}^{\oplus 3}_{S_{1}\times Y}\to \mathcal{E}'\to 
(\pi_{1}\times id)^{*}\mathcal{P}_{d} \to 0,\\
&0\to \mathcal{O}^{\oplus 4}_{S_{2}\times Y}\to \mathcal{F}'\to 
(\pi_{2}\times id)^{*}\mathcal{P}_{d}^{2} \to 0.
\end{align*}
Then $S=S_{1}\times_{J}S_{2},\: \mathcal{E}=
(\rho_{1}\times id)^{*}\mathcal{E}'\otimes 
p_{Y}^{*}\mathcal{O}_{Y}(-5Ny_{0})$ and 
$\mathcal{F}=(\rho_{2}\times id)^{*}\mathcal{F}'\otimes 
p_{Y}^{*}\mathcal{O}_{Y}(-8Ny_{0})$.
\end{proof}
The following statement is proved in \cite{AB} Lemma~10.1, when $g\geq 2$, 
and in \cite{Te} Lemma~2.3, when $g=1$. 
\begin{pro}\label{s2.26}
Let $F_{1},\ldots,F_{k}$ be semistable locally free sheaves on a 
smooth, projective curve. Then $F_{1}\otimes F_{2}\otimes \cdots 
\otimes F_{k}$ is semistable of slope $\mu(F_{1})+\mu(F_{2})+\cdots +\mu(F_{k})$. 
\end{pro}
\begin{cor}\label{s2.28}
Let $E$ and $F$ be semistable locally free 
sheaves on a smooth, projective curve $Y$. Then
$S^aE\otimes \bigwedge^bF$
is semistable of slope
\(
\mu = a\mu(E) + b\mu(F).
\)
\end{cor}
\begin{proof}
$E^{\otimes a}\otimes F^{\otimes b}$ is a semistable locally free sheaf and $S^aE\otimes \bigwedge^bF$ is one of its 
direct summands. So it is semistable of the same slope (see e.g. 
\cite{Oe} p.32).
\end{proof}
Using Proposition~\ref{s1.7} one obtains the following existence 
result for smooth covers of degree $d,\: 3\leq d\leq 5$.
\begin{pro}\label{s2.28a}
Let $d=3,\, 4$ or $5$. Let $Y$ be a smooth projective curve of genus 
$g\geq 1$. If $d=3$ let $E$ be a locally free sheaf of rank 2, 
which is stable in 
the case $g\geq 2$, and regular polystable in the case $g=1$. Let 
$R_{3}=S^3E\otimes (\det E)^{-1}$. If $d=4$ let $E$ and $F$ be locally 
free sheaves of ranks 3 and 2 respectively, 
which satisfy 
$\det F\cong \det E$, 
and
 which are 
stable if $g\geq 2$, and regular polystable if $g=1$. Let $\deg E=e$ and 
let $R_{4}=\Check{F}\otimes S^2E$. If $d=5$, let 
$E$ and $F$ be locally free sheaves of ranks 4 and 5 
 respectively, 
which satisfy 
$\det F\cong (\det E)^{2}$, 
and
which are stable if $g\geq 2$, and regular 
polystable if $g=1$. Let $\deg E=e$ and let $R_{5}=\bigwedge^2F\otimes 
E\otimes (\det E)^{-1}$. Suppose
\begin{equation}\label{es2.29}
e >
\begin{cases}
4(g-1)+4          & \text{if $d=3$}\\
12(g-1)+6         & \text{if $d=4$}\\
40(g-1)+20        & \text{if $d=5$}.
\end{cases}
\end{equation}
Then every $\eta$ in a Zariski open, dense subset of $H^0(Y,R_{d})$ is 
of the right codimension at every $y\in Y$ (see Section~\ref{s5}) and the 
associated finite covering $\pi_{\eta}:X_{\eta}\to Y$ of degree $d$ has
 smooth, irreducible $X_{\eta}$. 
\end{pro}
\begin{proof}
According to Corollary~\ref{s2.28} the locally free sheaf $R_{d}$ is 
semistable for every $d=3,\, 4$ or $5$. The inequalities 
\eqref{es2.29} are equivalent to $\mu(R_{d})>2g-1$. Therefore $R_{d}$ 
is globally generated (cf. \cite{Oe} p.39). 
Using the Bertini type theorems of \cite{CE} and \cite{C2} 
(cf. Proposition~\ref{s1.7}), 
there is a nonempty Zariski open subset in $H^0(Y,R_{d})$ such that 
each of its elements $\eta$ is of the right codimension at every 
$y\in Y$ and the corresponding cover $X_{\eta}\subset \mathbf{P}(E)$ 
is smooth. One has 
$(\pi_{\eta})_{*}\mathcal{O}_{X_{\eta}}\cong \mathcal{O}_{Y}\oplus 
\Check{E}$. The locally free sheaf $\Check{E}$ is stable, if $g\geq 
2$, and regular polystable if $g=1$, and  furthermore it is of negative 
degree. Hence $h^{0}(Y,\Check{E})=0$. It follows that  $X_{\eta}$ is 
connected, therefore irreducible.
\end{proof}
\begin{block}\label{s2.30}
Let $e\in \mathbb{Z}$. Let $d=3,\, 4$ or $5$. Let $S$ be the variety 
defined in Lemma~\ref{s2.20}, Lemma~\ref{s2.22} or Lemma~\ref{s2.24} 
respectively. Let $\mathcal{G}$ be $S^3\mathcal{E}\otimes (\det 
\mathcal{E})^{-1}$, $\Check{\mathcal{F}}\otimes S^2\mathcal{E}$ and 
$\bigwedge^2\mathcal{F}\otimes \mathcal{E}\otimes (\det 
\mathcal{E})^{-1}$ respectively. Suppose 
\begin{equation*}
e >
\begin{cases}
4(g-1)          & \text{if $d=3$}\\
12(g-1)         & \text{if $d=4$}\\
40(g-1)        & \text{if $d=5$}.
\end{cases}
\end{equation*}
If $s\in S$ let $\mathcal{G}_{s}=\mathcal{G}\otimes 
\mathcal{O}_{\{s\}\times Y}$. The inequalities are equivalent to 
$\mu(\mathcal{G}_{s})>2g-1$ for every $s\in S$. 
Since 
$\mathcal{G}_{s}$ is semistable (see Corollary~\ref{s2.28})
one has 
$h^1(Y,\mathcal{G}_{s})=0$ (see e.g. \cite{Oe} p.39). 
 The slope $\mu(\mathcal{G}_{s})$ 
takes values $\frac{e}{2},\, \frac{e}{6}$ and $\frac{e}{20}$ when 
$d=3,\, 4$ or $5$ respectively. Applying the Riemann-Roch theorem we 
obtain
\begin{equation}\label{es2.31}
h^0(Y,\mathcal{G}_{s}) =
\begin{cases}
2e +4(1-g)          & \text{if $d=3$}\\
2e+ 12(1-g)         & \text{if $d=4$}\\
2e + 40(1-g)        & \text{if $d=5$}.
\end{cases}
\end{equation}
Let $q:S\times Y\to S$ be the projection map. By Grauert's theorem 
$q_{*}\mathcal{G}$ is a locally free sheaf. We denote by 
$f:\mathbb{H}\to S$ the associated vector bundle with fibers 
$H^0(Y,\mathcal{G}_{s}),\, s\in S$. The subset $N\subset \mathbb{H}$, 
consisting of those $\eta$ which have the right codimension at every 
$y\in Y$, and furthermore have the property that 
$X_{\eta}\subset \mathbf{P}(\mathcal{E})_{f(\eta)}$
is smooth and irreducible, is Zariski open (possibly empty) in 
$\mathbb{H}$ (see Lemma~\ref{s5.71}). 
\end{block}
\begin{block}\label{s2.32}
Let us assume now that $e$ satisfies the stronger inequality 
\eqref{es2.29}. Using Proposition~\ref{s2.28a} we obtain that the 
subset $N\subset \mathbb{H}$ is nonempty and therefore dense in 
$\mathbb{H}$. Furthermore the map $f|_{N}:N\to S$ is surjective. Let 
$\mathcal{E}'$ (when $d=3$) and $\mathcal{E}',\mathcal{F}'$ (when 
$d=4$ or $5$) be the inverse images of $\mathcal{E}$ and $\mathcal{F}$ 
with respect to $f|_{N}\times id : N\times Y\to S\times Y$. 
 Let  $\mathbb{G}$ be the vector 
bundle over $S\times Y$ associated with the locally free sheaf $\mathcal{G}$. 
We apply the construction of  \S~\ref{s5.70} with $q:\mathcal{Y}\to \mathcal{Z}$ equal to the projection $S\times Y\to Y$. Here $p_1:\mathcal{Y}\times_{\mathcal{Z}}\mathbb{H}\to \mathcal{Y}$ may be identified with $f\times id:\mathbb{H}\times Y\to S\times Y$ and the fibers of the vector bundle $\mathbb{G}\to S\times Y$ are $\mathbb{G}_{(s,y)}=\mathcal{G}_s\otimes k(y)$. The diagram \eqref{es5.70} becomes
\[
\xymatrix{
              & \mathbb{G}\ar[d]\\
\mathbb{H}\times Y\ar[ru]^-{\Psi'}\ar[r]^-{f\times id}& S\times Y
}
\]
where $\Psi'(\eta,y)=\eta(y)\in \mathcal{G}_{f(\eta)}\otimes k(y)$. This yields a tautological 
section $\Psi \in H^0(\mathbb{H}\times Y,(f\times 
id)^{*}\mathcal{G})$. It has the property that for every $\eta \in 
\mathbb{H}$ one has $\Psi|_{\{\eta\}\times Y}=\eta \in 
H^0(Y,\mathcal{G}_{f(\eta)})$. 
For $d=3,4$
or $5$ the inverse image $(f|_{N}\times id)^{*}\mathcal{G}$ equals
respectively  $S^3\mathcal{E}'\otimes 
(\det 
\mathcal{E}')^{-1}$, 
$\Check{\mathcal{F}}'\otimes S^2\mathcal{E}'$ or 
$\bigwedge^2\mathcal{F}'\otimes \mathcal{E}'\otimes (\det 
\mathcal{E}')^{-1}$. The section
 $\Psi|_{N\times Y}$ has the
right codimension at every $(\eta,y)\in N\times Y$, 
so one obtains a Gorenstein covering $p:\mathcal{X}\to N\times Y$ of 
degree $d$ (see Section~\ref{s5}). For every $\eta \in N$   the restriction 
$p_{\eta}:\mathcal{X}_{\eta}\to \{\eta\}\times Y$ is equivalent to 
$\pi_{\eta}:X_{\eta}\to Y$ by Lemma~\ref{s5.71} and its Tschirnhausen module
is isomorphic to $\Check{\mathcal{E}}_s$, where $s=f(\eta)$.
The morphism $\pi_{N}\circ p:\mathcal{X}\to N$ is proper, smooth
with irreducible fibers
of dimension 1.
\par
Let $\mathcal{B}\subset N\times Y$ be the discriminant scheme of 
$p:\mathcal{X}\to N\times Y$. For each $\eta \in N$ the covering 
$\mathcal{X}_{\eta}\to \{\eta\}\times Y$ is ramified and the degree of 
the discriminant divisor is $n=2e$ (see e.g. \cite{K1} Lemma~2.3). 
The scheme $\mathcal{B}$ has pure 
codimension 1 in $N\times Y$ (\cite{AK} 
Theorem~VI.6.8). Since $N\times Y$ is smooth, $\mathcal{B}$ is a 
relative Cartier divisor in $N\times Y$. By the universal properties 
of $Y^{(n)}$ one obtains canonical 
morphism $\beta: N\to Y^{(n)}$ (cf \cite{ACGH} Ch.IV, Lemma~2.1). 
\end{block}
\begin{block}\label{s2.34}
We now want to focus on coverings with fixed determinants of the 
Tschirnhausen modules. Let $A\in Pic^{e}Y$.  For $d=3,\, 4$ or $5$  
let $S_{A}$ be the variety defined in Lemma~\ref{s2.20}, 
Lemma~\ref{s2.22} or Lemma~\ref{s2.24} respectively. If $\pi:X\to Y$ 
is a covering with $\pi_{*}\mathcal{O}_{X}\cong \mathcal{O}_{Y}\oplus 
\Check{E}$ and if $B$ is the discriminant divisor, then $(\det 
E)^{2}\cong \mathcal{O}_{Y}(B)$ (see \cite{Ha} Ch.IV, Ex.2.8 or 
\cite{AK} p.124).

Let the assumptions and notation be as in \S~\ref{s2.32}. If $\eta \in 
N$, let $X_{\eta}\to Y$ be the associated covering. Multiplying $\eta$ 
by a nonzero constant one obtains the same covering. So, $N$ is 
$\mathbb{C}^{*}$-invariant. Denoting by $N_{A}$ the preimage of 
$S_{A}$ one obtains morphisms
\[
\xymatrix{
\mathbb{P}N_{A}\ar[r]^{f}\ar[d]_{\beta}&S_A\\
|A^{2}|
}
\]
where $\beta(\langle \eta \rangle)=Discr(X_{\eta}\to Y)$. 
\end{block}
\begin{lem}\label{s2.35}
Let $d=3,\, 4$ or 5. Suppose $e$ satisfies the inequality 
\eqref{es2.29}. Let $A\in Pic^{e}Y$. Then $\beta: \mathbb{P}N_{A}\to 
|A^{2}|$ is a dominant morphism. 
\end{lem}
\begin{proof}
We want to prove that for every $\langle \eta \rangle\in 
\mathbb{P}N_{A}$ the fiber of $\beta:\mathbb{P}N_A\to |A^{2}|$ 
containing $\langle \eta \rangle$ has dimension equal to $\dim 
\mathbb{P}N_A-\dim |A^{2}|$ at $\langle \eta \rangle$. This is a local 
property, so it suffices to verify it for the complex analytic 
dimension of the intersection of the fiber through $\langle \eta 
\rangle$ with some complex neighborhood. Let $\langle \eta \rangle\in 
\mathbb{P}N_A$ and let $s=f(\langle \eta 
\rangle)$. Let $E$ (when $d=3$), and let $E,F$ (when $d=4$ or 5)       
be the vector bundles over $Y$ associated with $\mathcal{E}_{s}$ and 
$\mathcal{E}_{s}, \mathcal{F}_{s}$ respectively. Here 
$\mathcal{E},\mathcal{F}$ are the locally free sheaves over $S\times 
Y$ defined in Lemma~\ref{s2.20}, Lemma~\ref{s2.22} and 
Lemma~\ref{s2.24}. 
\par
Suppose first $g(Y)\geq 2$. Then $E$ and $F$ are stable bundles by 
construction. In particular they are simple. Using \cite{NS} 
Theorem~1 and the openness of the stability condition (see \cite{NS} 
Theorem~3) 
there is 
 a complex manifold $U_{A}$ and a holomorphic 
family of stable vector bundles $\{E''_{u}\}_{u\in U_{A}}$, 
respectively a pair of families $\{E''_{u}\}_{u\in U_{A}}$, 
$\{F''_{u}\}_{u\in U_{A}}$, 
with the property that $\det E_{u}''\cong A$ (when $d=3$), 
$\det E_{u}''\cong A\cong \det F_{u}''$ (when $d=4$), 
$\det E_{u}''\cong A$, $\det F_{u}''\cong A^{2}$ (when $d=5)$, the 
families are locally universal for deformations of $E$, respectively 
deformations of the pair $(E,F)$ with fixed determinants, and moreover
one may choose $U_{A}$ in such a way that different points of $U_{A}$
correspond to non-isomorphic bundles (when $d=3$), resp. pairs of 
bundles (when $d=4$ or $5$). So, there exists a complex neighborhood
$V\subset S_{A}$ of $s$ and a holomorphic map $g:V\to U_{A}$ such that
the families $\{E_{s}\}_{s\in V}$, resp. $\{E_{s}\}_{s\in V}, 
\{F_{s}\}_{s\in V}$ are pull-backs by $g$ of $\{E''_{u}\}_{u\in U_{A}}$
resp. $\{E''_{u}\}_{u\in U_{A}}$, $\{F''_{u}\}_{u\in U_{A}}$.
 Statement~(v) 
of Lemma~\ref{s2.20}, Lemma~\ref{s2.22} and 
Lemma~\ref{s2.24}
implies that $g$ has 
surjective differential at every point of $V$ in each of the corresponding 
cases.
 Replacing $U_{A}$ and $V$ by 
smaller open subsets we may thus assume that $g:V\to U_{A}$ is 
surjective, all fibers are smooth, equidimensional, and of dimension 
$\dim S_{A}-\dim U_{A}$. Let $R_{d}=S^3E\otimes (\det E)^{-1}$, 
$\Check{F}\otimes S^2E$ or $\bigwedge^2F\otimes E\otimes (\det E)^{-1}$ 
when $d=3,\, 4$ or $5$ respectively. Let $z=g(s)$ and let 
$V_{z}=g^{-1}(z)$. The families $\{E_{s}\}_{s\in V}$, $\{F_{s}\}_{s\in 
V}$ become trivial when restricted to $V_{z}$. Therefore 
$\mathbb{H}|_{V_{z}}\cong H^0(Y,R_{d})\times V_{z}$. This implies that the 
restriction of $\beta:\mathbb{P}N_A\to |A^{2}|$ to 
$f^{-1}(V_{z})$ 
is given by $(\langle \mu\rangle,v)\mapsto 
Discr(X_{\mu}\to Y)$, where $\langle \mu\rangle\in 
\mathbb{P}H^0(Y,R_{d})$, $v\in V_{z}$. Hence $\dim S_{A}-\dim U_{A}$ is a 
lower bound of the dimensions of the non-empty fibers of 
$\beta:\mathbb{P}N_A\to |A^{2}|$ at every $\langle \eta\rangle \in 
\mathbb{P}N_A$.
\par
Let $B=\beta(\langle \eta \rangle)$. Let us see for which $\eta'\in N$ 
it 
is possible that $B=\beta(\langle \eta' \rangle)$. 
Let $n=2e$. If $B\in Y^{(n)}$, 
there is a finite number (modulo equivalence) of coverings $\pi:X\to 
Y$ of degree $d$, with smooth irreducible $X$, whose discriminant 
equals $B$. Indeed, the larger set of classes of equivalence of coverings, 
whose discriminant locus is contained in $D=Supp(B)$ is in bijective correspondence 
with the homomorphisms $m:\pi_1(Y\setminus D,y_0)\to S_d$ modulo conjugation in $S_d$, $y_0\in Y\setminus D$ (see \cite{Fu} Proposition~1.2).
By construction the 
cover $X_{\mu}$ is smooth for every $\mu \in N$. 
So, if $\beta(\langle \eta'\rangle)=B$, then the equivalence class 
$[X_{\eta'}\to Y]$ may vary among a finite number of choices.
Two equivalent coverings yield isomorphic bundles. Indeed, given
$\pi:X\to Y$ one has $E^{\vee}\cong 
\pi_*\mathcal{O}_{X}/\mathcal{O}_{Y}$ and in case $d=4$ or $5$ one has
$F\cong Ker(S^{2}E\to \pi_{*}\omega_{X/Y}^{\otimes 2})$ (see 
\cite{CE} Theorem~4.4 and \cite{C2} Theorem~3.8). It follows that if 
$\beta(\langle \eta'\rangle)=B$, then $g\circ f(\langle \eta'\rangle)$
may vary among a finite number of points of $U_{A}$. Replacing $V$ and 
$U_{A}$ by smaller open sets we obtain that if 
$\langle \eta'\rangle\in f^{-1}(V)$ and 
$\beta(\langle \eta'\rangle)=B$, then 
$g\circ f(\langle \eta'\rangle)=z=g\circ f(\langle \eta\rangle)$. As 
we saw above the fiber of the composition 
$f^{-1}(V)\overset{f}{\lto}V\overset{g}{\lto}U_{A}$ over $z\in U_{A}$ 
is isomorphic to $W\times V_{z}$, where $W\subset \mathbb{P}H^0(Y,R_{d})$ 
is a Zariski open dense subset which parameterizes coverings as in 
\S~\ref{s2.30}. 
The bundles $E=E_{z}$, $F=F_{z}$ are stable, so $Aut(E)\cong 
\mathbb{C}^{*}\cong Aut(F)$. Using Proposition~\ref{s5.89b} we obtain 
that different elements of $W$ yield non-equivalent coverings. Hence 
the fiber of $\beta:f^{-1}(V)\to |A^{2}|$ over $B$ is contained in 
$f^{-1}(V_{z})$ and consists of points $(\langle \mu\rangle,v)$, where
$v\in V_{z}$ and $\langle \mu\rangle$ may vary among a finite number 
of points of $W\subset \mathbb{P}H^0(Y,R_{d})$. 
This implies that the dimension of the fiber of $\beta: \mathbb{P}N_A\to 
|A^{2}|$ at the point $\langle \eta\rangle$ is less or equal to $\dim 
S_{A}-\dim U_{A}$. Taking into account the previous inequality we 
conclude that every nonempty fiber of $\beta: \mathbb{P}N_A\to |A^{2}|$ 
is equidimensional of dimension $\dim S_{A}-\dim U_{A}$. 
\par
Using the dimension formula of \cite{NS} Theorem~1\, (iv) one 
obtains
\[
\dim U_{A} =
\begin{cases}
3(g-1)     &\text{if $d=3$}\\
11(g-1)    &\text{if $d=4$}\\
39(g-1)    &\text{if $d=5$}. 
\end{cases}
\]
The vector bundle $:\mathbb{H}\to S$ has rank given by \eqref{es2.31}. 
Taking into account that $n=2e$ we obtain
\[
\dim \left[\text{fiber$(\mathbb{P}N_A\to S_{A})$}\right] + \dim U_{A} = 
n-g.
\]The inequality \eqref{es2.29} implies $\deg A^{2}>2g-2$, thus $\dim 
|A^{2}|=n-g$. The calculation of the dimensions of the fibers of 
$\beta:\mathbb{P}N_A\to |A^{2}|$, made above, yields
\[
\begin{split}
\dim \overline{Im(\beta)} 
    &= \dim \mathbb{P}N_A - (\dim S_{A} - \dim U_{A})\\
    &= (\dim \mathbb{P}N_A - \dim S_{A}) + \dim U_{A}\\
    &= n-g.
\end{split}
\]
Therefore $\beta: \mathbb{P}N_A\to |A^{2}|$ is a dominant morphism and 
every nonempty fiber of $\beta$ has dimension 
$\dim \mathbb{P}N_A - \dim 
|A^{2}|$. 
\par
Suppose now $g(Y)=1$. Then $E$ and $F$ are regular polystable bundles. 
Using Proposition~\ref{s1.18aa} and the arguments of the case $g\geq 
2$ we obtain a complex manifold $U_{A}$, a complex neighborhood $V$ of 
$s$ in $S_{A}$ and a holomorphic map $g:V\to U_{A}$ with the same 
properties as in the case $g\geq 2$. Let $z=g(s)$ and let 
$V_{z}=g^{-1}(z)$. Again the restriction of 
$\beta: \mathbb{P}N_A\to |A^{2}|$ 
to $f^{-1}(V_{z})$ is given by $(\langle 
\mu\rangle,v)\mapsto Discr(X_{\mu}\to Y)$, where $\mu \in H^0(Y,R_{d}), 
v\in V_{z}$. The difference with the case $g\geq 2$ is that here 
$Aut(E)/\mathbb{C}^{*}$ (when $d=3$), and 
$Aut(E)/\mathbb{C}^{*}\times Aut(F)/\mathbb{C}^{*}$ (when $d=4$ or 
5) might be nontrivial. 
According to Proposition~\ref{s5.89b} these 
 groups act with finite stabilizers on
$W\subset \mathbb{P}H^0(Y,R_{d})$, the elements of one orbit yield
equivalent coverings, while the elements of different orbits yield
non-equivalent coverings. The same argument as in the case $g(Y)\geq 2$
shows that the nonempty fibers of $\beta: \mathbb{P}N_A\to |A^{2}|$
are equidimensional of dimension
\begin{align*}
\dim S_{A}-\dim U_{A} + \dim Aut(E)/\mathbb{C}^{*}     &\quad\text{(if $d=3$)}\\
\dim S_{A}-\dim U_{A} + \dim Aut(E)/\mathbb{C}^{*} + \dim Aut(F)/\mathbb{C}^{*}
&\quad \text{(if $d=4$ or $5$).}
\end{align*}
Let $G=G_{1}\oplus \cdots \oplus G_{h}$ be the decomposition of a regular 
polystable bundle in a direct sum of stable bundles. Then $Aut(G)\cong 
(\mathbb{C}^{*})^{h}$, so $\dim Aut(G)/\mathbb{C}^{*}=h-1$. If $d=3$ let 
$h=g.c.d.(2,e)$. If $d=4$ let $h'=g.c.d.(3,e)$ and let 
$h''=g.c.d.(2,e)$. If $d=5$ let $h'=g.c.d.(4,e)$ and let 
$h''=g.c.d.(5,2e)$. 
We obtain that the dimension of the nonempty fibers of 
$\beta:\mathbb{P}N_A\to |A^{2}|$ equals
\begin{align*}
\dim S_{A}-\dim U_{A} +h -1          &\quad \text{(if $d=3$)}\\
\dim S_{A}-\dim U_{A} +h' -1+h''-1   &\quad \text{(if $d=4$ or $5$.)}
\end{align*}
The dimension of $U_{A}$ equals $h^1(Y,ad'E)=h-1$ if $d=3$ and 
$h^1(Y,ad'E)+h^1(Y,ad'F)=h'-1+h''-1$ if $d=4$ or $5$. According to 
\eqref{es2.31} one has $\dim \mathbb{P}N_A = n-1 +\dim S_{A}$, so we 
obtain
\begin{align*}
\dim \overline{Im(\beta)} &{} \\
=&\dim \mathbb{P}N_A - \dim \left[ \text{fiber($\mathbb{P}N_A \to 
|A^{2}|$)}\right]&{}\\
=&n-1+\dim S_{A}-(\dim S_{A}-\dim U_{A} +h-1)&\text{(if $d=3$)}\\
=&n-1+\dim S_{A}-(\dim S_{A}-\dim U_{A} +h'-1+h''-1) &\text{(if $d=4$ or 
5)}\\
=&n-1. &{}
\end{align*}
Therefore $\beta : \mathbb{P}N_A\to |A^{2}|$ is a dominant morphism and 
every nonempty fiber of $\beta$ has dimension 
$\dim \mathbb{P}N_A - \dim 
|A^{2}|$. 
\end{proof} 
\begin{lem}\label{s2.40}
Let $\pi:X\to Y$ be a covering of degree $d$ of smooth, irreducible, 
projective curves. Let $g(Y)\geq 1$. Suppose $\pi$ is simply branched 
in $n>0$ points and furthermore suppose the Tschirnhausen module 
$\Check{E}=\pi_{*}\mathcal{O}_{X}/\mathcal{O}_{Y}$ is semistable. Then 
the monodromy group of the covering is  the symmetric 
group $S_{d}$. 
\end{lem}
\begin{proof}
Suppose, by way of contradiction, that the monodromy group is a proper 
subgroup of $S_{d}$. According to \cite{BeE} Lemma~2.4 the map $\pi$ 
can be factored as 
$X\overset{\pi_1}{\to}\tilde{Y}\overset{\pi_2}{\to}Y$ with $\deg 
\pi_{1}>1$, $\deg \pi_{2}>1$. The covering $\pi_{2}:\tilde{Y}\to Y$ is 
\'{e}tale since $\pi$ is simply branched. One has an injective 
homomorphism of locally free sheaves \linebreak
$\pi_{1}^{*}:\pi_{2*}\mathcal{O}_{\tilde{Y}}\to 
\pi_{*}\mathcal{O}_{X}$, which is identity on $\mathcal{O}_{Y}$. One 
obtains an injective homomorphism 
\(
\pi_{2*}\mathcal{O}_{\tilde{Y}}/\mathcal{O}_{Y}\lto 
\pi_{*}\mathcal{O}_{X}/\mathcal{O}_{Y} = \Check{E}.
\)
This is absurd since 
$\deg(\pi_{2*}\mathcal{O}_{\tilde{Y}}/\mathcal{O}_{Y})=0$ and 
$\Check{E}$ is semistable of negative degree $-\frac{n}{2}$.
\end{proof}

\begin{lem}\label{s2.41}
Let the assumptions  be as in Lemma~\ref{s2.35}. 
The subset $M\subset N$ consisting of $\eta$ such that the covering 
$\pi_{\eta}:X_{\eta}\to Y$ has smooth, irreducible $X_{\eta}$, simple 
branching and monodromy group equal to $S_{d}$, is 
$\mathbb{C}^{*}$-invariant, Zariski open and dense in $N$. For every
$A\in Pic^{e}Y$ the intersection $M_{A}=M\cap N_{A}$ is Zariski 
open and dense in $N_{A}$. 
\end{lem}
\begin{proof}
By Lemma~\ref{s2.40} the subset $M\subset N$ consists of points $\eta$ 
such that $\pi_{\eta}:X_{\eta}\to Y$ is simply branched. This subset 
is open according to Lemma~\ref{s5.71} and is obviously 
$\mathbb{C}^{*}$-invariant. Let $A\in Pic^{e}Y$. The variety $N_{A}$ 
is irreducible since it is an open subset of a vector bundle over the 
irreducible variety $S_{A}$. By Lemma~\ref{s2.35} and Bertini's 
theorem applied to $|A^{2}|$ we conclude that 
$M_{A}=M\cap N_{A}\neq \emptyset$, hence it is dense in $N_{A}$. 
Therefore $M$ is dense in $N$ as well.
\end{proof}
\section{Unirationality of Hurwitz spaces}\label{s3}
\begin{block}\label{s3.51}
Let $\mathcal{H}_{d,n}(Y)$ be the Hurwitz space, which 
parameterizes the equivalence classes $\{ [\pi:X\to Y]\}$ of  degree $d$ 
coverings of $Y$ simply branched in $n$ points, with  smooth, irreducible $X$.
This space is a finite unramified cover of $Y^{(n)}\setminus \Delta$, where
$\Delta$ is the codimension one subvariety parameterizing divisors 
with multiplicities, the morphism $b: \mathcal{H}_{d,n}(Y)\to Y^{(n)}\setminus \Delta$ associates to $[\pi:X\to Y]$ its discriminant divisor.
A structure of a complex analytic manifold on  $\mathcal{H}_{d,n}(Y)$ such that $b$ is a finite \'{e}tale analytic map is given  in \cite{Fu} Sect.1. That $\mathcal{H}_{d,n}(Y)$
has a structure of algebraic variety, such that $b$ is a morphism, and moreover such a structure 
is unique up to isomorphisms compatible with $b$, follows from the 
generalized Riemann's Existence Theorem (\cite{SGA1} Expos\'{e} XII,  Th\'{e}or\`{e}me~5.1).
We denote by 
$\mathcal{H}_{d,n}^0(Y)$ the union of connected components of 
$\mathcal{H}_{d,n}(Y)$, which parameterize equivalence classes of 
coverings, whose monodromy group is $S_d$.
\end{block} 
\begin{lem}\label{s3.51a}
Let $X$ and $Y$ be smooth, irreducible, projective curves
and let $p:X\to Y$ be a covering of degree $d$ whose monodromy group 
is $S_{d}$. Then $Aut(X/Y)=\{1\}$.
\end{lem}
\begin{proof}
Let $D$ be the branch locus of $p$, let $Y_0=Y\setminus D$ and let $X_0=p^{-1}(Y_0)$. Let 
$y_0\in Y_0$, $p^{-1}(y_0)=\{x_1,\ldots,x_d\}$. Let $H=p_*\pi_1(X_0,x_1)\subset \pi_1(Y_0,y_0)$.
By \cite{Mas} Ch.~V Cor.~7.3 $Aut(X_0/Y_0)\cong N[H]/H$, where $N[H]$ is the normalizer of $H$ in $\pi_1(Y_0,y_0)$. The monodromy (right) action of $\pi_1(Y_0,y_0)$ on $\{x_1,\ldots,x_d\}$ determines a $\pi_1(Y_0,y_0)$-equivariant bijection $H\backslash \pi_1(Y_0,y_0) \to \{x_1,\ldots,x_d\}$. Hence the kernel of the monodromy homomorphism 
$m:\pi_1(Y_0,y_0)\to S_d$ is  $\cap gHg^{-1}\subset H$. The image of $H$ is the subgroup $G_1\subset S_d$ of permutations which fix $x_1$. The homomorphism $m$ induces an isomorphism $N[H]/H\cong N[G_1]/G_1$, where $N[G_1]$ is the normalizer of $G_1$ in $S_d$. Clearly $N[G_1]=G_1$, so $Aut(X/Y)\cong Aut(X_0/Y_0) = \{id\}$.
\end{proof}
\begin{block}\label{fine moduli}
The algebraic variety $\mathcal{H}^{0}_{d,n}(Y)$ is a fine moduli scheme for simple coverings of $Y$ with full monodromy group $S_d$. This is proved as follows. The existence of fine moduli scheme $H$ which is an \'{e}tale cover of  
$Y^{(n)}\setminus \Delta$ can be deduced from the results in \cite{We} (see Theorem~4 and also \cite{RW}).  The same argument as in \cite{Fu} Proposition~7.3 shows that there is an analytic isomorphism of $H$ with 
$\mathcal{H}^{0}_{d,n}(Y)$ compatible with the two covering morphisms to $Y^{(n)}\setminus \Delta$.
The uniqueness result of \cite{SGA1} Expos\'{e} XII,  Th\'{e}or\`{e}me~5.1 shows that this is an algebraic isomorphism.
\par
We mention the approach to the Hurwitz spaces, based on the Kontsevich moduli spaces, developed in \cite{HGS}. As the referee points out
one can work with the connected moduli scheme $H^{d,n}_{S_d}(Y,\Sigma, y_0)$ described in this paper.
\par
Let $\mathcal{X}\to 
\mathcal{H}^{0}_{d,n}(Y)\times Y$ be the universal  covering.  Let 
$\mathcal{E}$ be the locally free sheaf of rank $d-1$ dual to the 
Tschirnhausen module of the universal covering. Let $\mathcal{L}=
\det(\mathcal{E})$. For every $\xi \in \mathcal{H}^{0}_{d,n}(Y)$ one has 
that $\deg (\mathcal{L}|_{\{\xi\}\times Y}) = \deg(\mathcal{E}_{\xi}) = 
\frac{n}{2}$. Let $e=\frac{n}{2}$. By the universal property of 
$Pic^{e}Y$, there is a morphism $t:\mathcal{H}^{0}_{d,n}(Y)\to Pic^{e}Y$, 
such that $t(\xi)=\det(\mathcal{E}_{\xi})$. For every $A\in 
Pic^{e}(Y)$ we denote by $\mathcal{H}^{0}_{d,A}(Y)$ the reduced scheme 
$t^{-1}(A)$.  
\end{block}
\begin{block}\label{s3.52}
Let $d=3,\, 4$ or $5$. Suppose $e$ satisfies the inequality 
\eqref{es2.29}. We use the notation of \S~\ref{s2.32}. Let 
$M\subset N$ and $M_{A}\subset N_{A}$ be the open dense subsets 
defined in Lemma~\ref{s2.41}.
 Restricting the covering $\mathcal{X}\to N\times 
Y$ to $M\times Y$ and $M_{A}\times Y$, and using the universal 
property of the Hurwitz spaces, one obtains canonical morphisms 
$h: \mathbb{P}M\to \mathcal{H}_{d,n}^0(Y)$ and $h_{A}: \mathbb{P}M_{A}\to 
\mathcal{H}_{d,A}^0(Y)$.
\end{block}
\begin{proof}[Proof of Theorem~\ref{s0.3}]
We first prove the theorem for
general $A\in 
Pic^{e}Y$.
The hypothesis on $e$ implies that $n=2e\geq 2d$. According to 
\cite{BeE} Theorem~6.3, see also \cite{HGS}, the Hurwitz space $\mathcal{H}_{d,n}^0(Y)$ is 
connected. Therefore $\mathcal{H}_{d,n}^0(Y)$ is irreducible, 
since it is a smooth variety. The morphism $h: \mathbb{P}M\to 
\mathcal{H}_{d,n}^0(Y)$ fits into the following commutative diagram. 
\begin{equation}\label{es3.54}
\xymatrix{
\mathbb{P}M\ar[d]\ar[r]^-{h}&\mathcal{H}_{d,n}^0(Y)\ar[d]_{t}\ar[r]^-{b}
&Y^{(n)}\setminus \Delta\ar[d]\\
Pic^{e}Y\ar@{=}[r]&Pic^{e}Y\ar[r]^{sq}&Pic^{n}Y.
}
\end{equation}
Here 
the morphism $b$ is defined in \S~\ref{s3.51},
$sq(A)=A^{2}$, the morphism $t$ is defined in 
\S~\ref{fine moduli}, the left vertical morphism is defined similarly, and 
the right vertical one is $D\mapsto \mathcal{O}_{Y}(D)$. The 
composition $b\circ h$  equals the restriction of
$\beta:\mathbb{P}N\to Y^{(n)}$ to $\mathbb{P}M$ (see \S~\ref{s2.32}). Using 
Lemma~\ref{s2.35} and Lemma~\ref{s2.41} we conclude that this 
composition is dominant. This implies that $h: \mathbb{P}M\to 
\mathcal{H}_{d,n}^0(Y)$ is dominant as well, since $b$ is a finite 
\'{e}tale covering, and $\mathcal{H}_{d,n}^0(Y)$ is irreducible. 
Let $L\in Pic^{n}Y$. One has $\deg L = n = 2e \geq 2g+1$, so 
$\dim |L|=n-g$ and $Y^{(n)}\setminus \Delta\to Pic^{n}Y$ is a smooth 
morphism. One has
\begin{equation}\label{es3.54a}
b^{-1}(|L|\setminus \Delta) = \bigcup_{A\in Pic^{e}Y,\; A^{2}\cong L}
\mathcal{H}_{d,A}^0(Y).
\end{equation}
For every $A\in Pic^{e}Y$ with $A^{2}\cong L$ the set 
$\mathcal{H}_{d,A}^0(Y)$ is nonempty by Lemma~\ref{s2.41}. Therefore 
$\mathcal{H}_{d,A}^0(Y)$ is a smooth, equidimensional variety of 
dimension $n-g$.
Let $Z$ be the closure in $\mathcal{H}_{d,n}^0(Y)$ of 
the complement to the image of $\mathbb{P}M$. One has $\dim Z<n$, so 
there are two possibilities. Either $t|_{Z}:Z\to Pic^{e}Y$ is not 
dominant, or $t|_{Z}:Z\to Pic^{e}Y$ is dominant, but the general fiber 
has dimension $<n-g$. In either case
$h_{A}:\mathbb{P}M_A\to \mathcal{H}_{d,A}^0(Y)$ is dominant
 for every sufficiently general 
$A\in Pic^{e}Y$. Using the fact that $M_{A}$ is a 
Zariski open dense subset in a vector bundle over the rational variety 
$S_{A}$ we conclude that $\mathcal{H}_{d,A}^0(Y)$ is irreducible and 
unirational. This proves the theorem when $g(Y)\geq 1$ and $A$ is 
general. When $g(Y)= 
1$ one has $\mathcal{H}_{d,A}^0(Y) \cong \mathcal{H}_{d,A'}^0(Y)$ for 
every pair $A,A'\in Pic^{e}Y$ (cf. \cite{K1} Lemma~2.5). Therefore 
$\mathcal{H}_{d,A}^0(Y)$ is irreducible and unirational for every 
$A\in Pic^{e}Y$.
\par
It remains to prove Theorem~\ref{s0.3} when $g(Y)\geq 2$ and $A\in 
Pic^{e}Y$ is arbitrary.
 Let $L\in Pic^{n}Y$ be arbitrary. It suffices 
to prove that in the disjoint union \eqref{es3.54a} every 
$\mathcal{H}_{d,A}^0(Y)$ is connected. Indeed, then 
$\mathcal{H}_{d,A}^0(Y)$ is irreducible, the composition 
$\mathbb{P}M_{A}\to \mathcal{H}_{d,A}^0(Y)\overset{b}{\to}|A^{2}|$ is 
dominant by Lemma~\ref{s2.35} and Lemma~\ref{s2.41}. 
The morphism  $b$ is an \'{e}tale covering. Hence the 
morphism $h_{A}:\mathbb{P}M_{A}\to \mathcal{H}_{d,A}^0(Y)$ is dominant. 
Therefore $\mathcal{H}_{d,A}^0(Y)$ is unirational.
\par
The statement that the smooth variety $\mathcal{H}_{d,A}^0(Y)$ is 
connected for every $A\in Pic^{e}Y$ with $A^{2}\cong L$ is equivalent 
to the statement that $b^{-1}(|L|\setminus \Delta)$ has $2^{2g}$ 
connected components. This is clear from \eqref{es3.54a} since
$sq:Pic^eY\to Pic^nY$, $sq(A)=A^2$ is up to isomorphism the same map as the isogeny
$2\cdot id_J:J(Y)\to J(Y)$ of degree $2^{2g}$.
We know the above property holds for every $L'$ in a Zariski 
open dense subset $U\subset Pic^{n}Y$. Let $D\in |L|\setminus \Delta$.
Let $G$ be the image of $\pi_{1}(|L|\setminus \Delta, D)\to
\pi_{1}(Y^{(n)}\setminus \Delta, D)$. The monodromy action of 
$\pi_{1}(|L|\setminus \Delta, D)$ on $b^{-1}(D)$ factors through the 
action of $G$. So, the number of connected components of 
$b^{-1}(|L|\setminus \Delta)$ equals the number of orbits of $G$
with respect to the action on $b^{-1}(D)$. According to a result due 
to Dolgachev and Libgober (cf. Proposition~\ref{s4.60}) one has an exact 
sequence
\begin{equation*}
0\to \pi_{1}(|L|\setminus \Delta,D)\to 
\pi_{1}(Y^{(n)}\setminus \Delta,D)\to 
\pi_{1}(Pic^{n}Y,L)\to 1.
\end{equation*}
So, $G=Ker(\pi_{1}(Y^{(n)}\setminus \Delta,D)\to 
\pi_{1}(Pic^{n}Y,L))$. Let $L'\in U\subset Pic^{n}Y$ and let $D'\in 
|L'|$. We know that
$G'=Ker(\pi_{1}(Y^{(n)}\setminus \Delta,D')\to 
\pi_{1}(Pic^{n}Y,L'))$ has $2^{2g}$ orbits with respect to the action 
of $G'$ on $b^{-1}(D')$. Connecting $D$ and $D'$ by a path in 
$Y^{(n)}\setminus \Delta$ we have a commutative diagram
\begin{equation*}
\xymatrix{
\pi_{1}(Y^{(n)}\setminus \Delta,D)
\ar[r]\ar[d]&\pi_{1}(Pic^{n}Y,L)\ar[d]^{\cong}\\
\pi_{1}(Y^{(n)}\setminus \Delta,D')\ar[r]&\pi_{1}(Pic^{n}Y,L')
}
\end{equation*}
where the left vertical isomorphism is compatible with the actions
on $b^{-1}(D)$ and $b^{-1}(D')$ respectively. We conclude that 
$G$ and $G'$ have the same number of orbits. Therefore 
$b^{-1}(|L|\setminus \Delta)$ has $2^{2g}$ connected components.
Theorem~\ref{s0.3} is proved.
\end{proof}
\begin{proof}[Proof of Theorem~\ref{s0.100}]
Let us first consider the case $g(Y)\geq 2$.  We 
recall that when $(r,e)=1$ there is a smooth, irreducible, projective 
variety $U(r,e)$ which is a coarse moduli space for stable vector 
bundles of rank $r$ and degree $e$. Furthermore there is a Poincar\'{e} locally 
free sheaf $\mathcal{P}(r,e)$ on
$U(r,e)\times Y$ (see \cite{T2} or \cite{N2} Ch.5 \S~5).  Let $A\in Pic^{e}Y$, let $SU(r,A)$ 
be the closed 
subvariety which parameterizes the stable vector bundles with 
determinant isomorphic to $A$ and let $\mathcal{P}(r,A)$ be the 
restriction of $\mathcal{P}(r,d)$ on $SU(r,A)$. 
The arguments of \S~\ref{s2.30}, \S~\ref{s2.32} and 
\S~\ref{s2.34}  may be applied to new families of vector bundles on $Y$.
Namely: let $S_{A}=SU(2,A)$, 
$\mathcal{E}=\mathcal{P}(2,A)$ when $d=3$; let $S_{A}=SU(3,A)\times SU(2,A)$, 
\; $\mathcal{E}$ and $\mathcal{F}$ be respectively equal to the inverse images of 
$\mathcal{P}(3,A)$ and $\mathcal{P}(2,A)$  when $d=4$; 
$S_{A}=SU(4,A)\times SU(5,A^{2})$, \; 
$\mathcal{E}$ and $\mathcal{F}$ be respectively equal to the inverse images of 
$\mathcal{P}(4,A)$ and $\mathcal{P}(5,A^{2})$  when $d=5$.
For $d=3,4$ or $5$  let 
 $\mathcal{G}$ be $S^3\mathcal{E}\otimes (\det 
\mathcal{E})^{-1}$, $\Check{\mathcal{F}}\otimes S^2\mathcal{E}$ and
$\bigwedge^2\mathcal{F}\otimes \mathcal{E}\otimes (\det 
\mathcal{E})^{-1}$ respectively.
One repeats all the arguments after \S~\ref{s2.32} concerning 
$\mathbb{P}N_{A}, \mathbb{P}M_{A}$ etc. with the new 
definition of $S_{A}$. 
We notice that the compatibility conditions $\det \mathcal{E}\cong \det \mathcal{F}$
(when $d=4$) and $(\det \mathcal{E})^{ 2}\cong \det \mathcal{F}$ (when $d=5$)
might not be satisfied, but for every $s\in S_{A}$ one has 
$\det \mathcal{E}_{s}\cong A \cong \det \mathcal{F}_{s}$ and 
$(\det \mathcal{E}_{s})^{ 2}\cong A^{2} \cong \det \mathcal{F}_{s}$ respectively, and this suffices for applying the arguments after \S~\ref{s2.32}. 
In particular the proof of Lemma~\ref{s2.35}
simplifies. One may let $V=S_{A}=U_{A}$ and $g$ be the identity map.
One obtains a morphism $h_{A}:\mathbb{P}M_{A}\to \mathcal{H}_{d,A}^0(Y)$
whose composition with the finite \'{e}tale covering 
$b_{A}:\mathcal{H}_{d,A}^0(Y)\to |A^{2}|$ is 
dominant. 
According to Theorem~\ref{s0.3} the variety $\mathcal{H}_{d,A}^0(Y)$ is irreducible.
Therefore $h_{A}:\mathbb{P}M_{A}\to \mathcal{H}_{d,A}^0(Y)$ is dominant.
We claim that this morphism is moreover injective. Indeed, given 
a $[\pi:X\to Y]\in \mathcal{H}_{d,A}^0(Y)$ in the image of $h_{A}$, the 
locally free sheaves $E$ (when $d=3$) and $E$, $F$ (when $d=4$ or $5$) 
are uniquely determined, up to isomorphism, from the covering since
$E^{\vee}\cong \pi_*\mathcal{O}_{X}/\mathcal{O}_{Y}$ and 
$F\cong Ker(S^{2}E\to \pi_{*}\omega_{X/Y}^{\otimes 2})$ (see 
\cite{CE} Theorem~3.4 and Theorem~4.4, and \cite{C2} Theorem~3.8). 
The fiber $W$ of 
$h_{A}:\mathbb{P}M_{A}\to S_{A}$ over $[E]$ (when $d=3$) or $([E],[F])$ 
(when $d=4$ or $5$) may be identified with the Zariski open dense 
subset of 
$\mathbb{P}H^0(Y,R_{d})$ consisting of $\langle \eta \rangle$ such 
that: $X_{\eta}$ is smooth and irreducible; $\pi_{\eta}:X_{\eta}\to Y$ 
is simply ramified; its monodromy group is $S_{d}$ (see 
Lemma~\ref{s5.71}, Proposition~\ref{s2.28a}, Lemma~\ref{s2.40} and 
Lemma~\ref{s2.41}). According to Lemma~\ref{s3.51a} for every 
$\langle \eta \rangle \in W$ one has $Aut(X_{\eta}/Y)=\{1\}$. Using 
Proposition~\ref{s5.89b} and the fact that $Aut(E)\cong 
\mathbb{C}^{*}$, $Aut(F)\cong \mathbb{C}^{*}$ we conclude that 
$h_{A}|_{W}:W\to \mathcal{H}_{d,A}^0(Y)$ is injective. This proves that
$h_{A}:\mathbb{P}M_{A}\to \mathcal{H}_{d,A}^0(Y)$ is injective, therefore 
$h_{A}$ is a birational isomorphism. The variety $S_{A}$, being a product 
of moduli spaces $SU(r,A)$ with $(r,\deg A)=1$, is rational according 
to the theorem of King and Schofield \cite{KiS} Theorem~1.2. We should notice that
one may use results due to Tyurin and Newstead when $d=3$ or $4$: \cite{T1} Theorem~11 in 
the case $d=3$ and \cite{N1} Proposition~2 in the case $d=4$.
The variety $\mathbb{P}M_{A}$ is a Zariski open subset of the 
projectivization of a vector bundle over $S_{A}$. Therefore 
$\mathbb{P}M_{A}$ and $\mathcal{H}_{d,A}^0(Y)$ are rational varieties.
\par
Suppose now $g(Y)=1$. If $(r,e)=1$ and $A\in Pic^{e}Y$ there is a 
unique, up to isomorphism, indecomposable (= stable) locally free 
sheaf of rank $r$ and determinant isomorphic to $A$ (\cite{At} p.434). 
 We may apply the same arguments as in the case $g(Y)\geq 
2$ taking for $SU(r,A)$ and $S_{A}$ varieties consisting of one point.
Theorem~\ref{s0.100} is proved.
\end{proof}
\appendix
\section{A result of Dolgachev and Libgober}
Let $Y$ be a smooth, irreducible, projective curve of genus $g\geq 
1$. Let $Y^{(n)}\to Pic^{n}Y$ be the canonical morphism 
$u(D)=\mathcal{O}_{Y}(D)$. Let $\Delta \subset Y^{(n)}$ be the 
codimension one subvariety which parameterizes divisors with 
multiplicities. We need the following statement from \cite{DL} p.9.
\begin{pro}[Dolgachev-Libgober]\label{s4.60}
Let $n\geq 2g+1$. Then the canonical morphism $u:Y^{(n)}\setminus 
\Delta \to Pic^{n}Y$ is a Serre fibration. For every 
$L\in Pic^{n}Y$ and every $D\in |L|$ one has an exact sequence
\begin{equation}\label{es4.60}
1\to \pi_{1}(|L|\setminus \Delta, D)\to 
\pi_{1}(Y^{(n)}\setminus \Delta, D) \to 
\pi_{1}(Pic^{n}Y, L)\to 1
\end{equation}
\end{pro}
A weaker result, \eqref{es4.60} under the assumption that $L$ is 
general enough, is  published in \cite{Sh}, where it is argued 
that the proof in \cite{DL} seems to be incomplete (see \cite{Sh} p.337). We 
need the exactness of \eqref{es4.60} for every $L$, so we include 
below
a detailed proof of Proposition~\ref{s4.60} along the sketch given
in \cite{DL} p.9.
\begin{lem}\label{s4.60a}
Let
\begin{equation*}
\xymatrix{
A\ar[rr]^{h}\ar[rd]_{p}&&
B\ar[dl]^{q}\\
&S
}
\end{equation*}
be a commutative diagram of continuous maps of topological spaces. 
Assume $p$ and $h$ are Serre fibrations and $h$ is surjective. Then 
$q$ is a Serre fibration.
\end{lem}
\begin{proof}
Let $I=[0,1]$ be the unit interval. Let $f:I^{n}\to B$ be a continuous
map and let $G:I^{n}\times I\to S$ be a homotopy such that 
$G|_{I^{n}\times \{0\}}=q\circ f$. Let $b=f(0,\ldots,0)$. There is a 
point $a\in A$ such that $h(a)=b$. Since $h$ is a Serre fibration 
there is a lifting $\tilde{f}:I^{n}\to A$ such that 
$h\circ \tilde{f}=f$. One has $p\circ \tilde{f}=q\circ f$. Since 
$p:A\to S$ is a Serre fibration there is a homotopy 
$\tilde{F}:I^{n+1}\to A$ which lifts $G:I^{n+1}\to S$. One lets 
$F=h\circ \tilde{F} :I^{n+1}\to B$. This homotopy lifts $G$.
\end{proof}
\begin{block}\label{s4.61}\emph{Families of dual varieties}
The classical arguments about dual varieties (see e.g. \cite{Lo} 
Lecture~1) are easily extended to families of smooth varieties. Let 
$\mathcal{X}$ and $S$ be smooth irreducible varieties and let 
$f:\mathcal{X}\to S$ be a 
proper, smooth morphism with irreducible fibers. Let $V\to S$ be a 
vector bundle of rank $N+1$ and let $\varphi :\mathcal{X}\to 
\mathbb{P}(V)$ be a $S$-morphism, which is a closed embedding. One considers 
the dual vector bundle $V^{*}$, the incidence variety 
$W\subset \mathbb{P}(V)\times_{S}\mathbb{P}(V^{*})$ and the two 
projections $p:W\to \mathbb{P}(V)$ and $q:W\to \mathbb{P}(V^{*})$.
One lets $W_{\mathcal{X}}=p^{-1}(\varphi(\mathcal{X}))$, 
$q_{\mathcal{X}}:=q|_{W_{\mathcal{X}}}$ and 
\begin{equation*}
N_{\mathcal{X}}= \{(\varphi(x),\alpha)\in W_{\mathcal{X}}
\mid \mathbb{H}_{\alpha}\supset d\varphi\, (T_{\mathcal{X}/S}(x))\}.
\end{equation*}
One proves as in \cite{Lo} Lecture~1 that 
$W_{\mathcal{X}}$ and $N_{\mathcal{X}}$ are smooth fibrations over $S$ 
of relative dimensions $\dim \mathcal{X}_{s}+rk(V)-2$ and 
$rk(V)-2$ respectively.

\smallskip \noindent
CLAIM: \emph{
The set of critical points of $q_{\mathcal{X}}:W_{\mathcal{X}}\to 
\mathbb{P}(V^{*})$  equals  $N_{\mathcal{X}}$ and  the set 
of critical values of $q_{\mathcal{X}}$ equals 
$\cup_{s\in S}\varphi_{s}(\mathcal{X}_{s})^{\vee}$. 
The map 
$q_{\mathcal{X}}:W_{\mathcal{X}}\to \mathbb{P}(V^{*})$ is locally trivial in the $C^{\infty}$-category over 
$\mathbb{P}(V^{*})\setminus \cup_{s\in S}
\varphi_{s}(\mathcal{X}_{s})^{\vee}$.}
\begin{proof}
Let $a\in \mathcal{X}_{s}$, $\alpha \in \mathbb{P}(V^{*}_{s})$ satisfy $a\in \mathbb{H}_{\alpha}$, i.e. $(a,\alpha)\in (W_{\mathcal{X}})_{s}$. One has a commutative diagram with exact rows
\[
\xymatrix{
0\ar[r]&T_{W_{\mathcal{X}}/S}(a,\alpha)\ar[d]^{dq_{\mathcal{X}}}\ar[r]
&T_{W_{\mathcal{X}}}(a,\alpha)\ar[d]^{dq_{\mathcal{X}}}\ar[r]
&T_{S}(s)\ar[d]^{=}\ar[r]&0\\
0\ar[r]&T_{\mathbb{P}(V^{*})/S}(\alpha)\ar[r]&T_{\mathbb{P}(V^{*})}(\alpha)\ar[r]&T_{S}(s)\ar[r]&0
}
\]
The middle vertical map is surjective if and only if the left vertical map is surjective.
Hence the first part of the claim follows from \cite{Lo} Lemma~1.4 and the second one from the Ehresmann fibration theorem (see \cite{Eh} p.154)
\end{proof}
\end{block}
\begin{block}\label{s4.62}\emph{Example.}
Let $C$ be a smooth, projective, irreducible curve of genus $g\geq 0$. 
Let $n\geq 2g+1$. Let $S=Pic^{n}C$, $\mathcal{X}=C\times S$ and 
$f:\mathcal{X}\to S$ be the second projection. Let $\mathcal{L}$ be a 
Poincar\'{e} invertible sheaf on $\mathcal{X}$ and let 
$E=f_{*}\mathcal{L}$ be the associated locally free sheaf of rank 
$n-g+1$. Let $V$ be the vector bundle associated with $E^{\vee}$ with 
fibers $H^0(C,\mathcal{L}_{s})^{*}$ and let 
$\varphi :C\times S\to \mathbb{P}(V)$ be the embedding 
determined by $\mathcal{L}$. Then $\mathbb{P}(V^{*})\to S$  is 
identified with the canonical morphism $u:C^{(n)}\to Pic^{n}C$ (see 
\cite{ACGH} p.309), $W_{\mathcal{X}}$ and $N_{\mathcal{X}}$ are 
identified with the following subsets of $C\times C^{(n)}$:
\begin{equation*}
W_{\mathcal{X}}=\{(x,D)\mid x\in Supp(D)\},\quad 
N_{\mathcal{X}} = \{(x,D)\mid \nu_{x}(D)\geq 2\}.
\end{equation*}
One concludes that the image of $W_{\mathcal{X}}\setminus 
N_{\mathcal{X}}$ in $C^{(n)}$ equals $C^{(n)}\setminus \Delta$. 
\end{block}
\begin{proof}[Proof of Proposition~\ref{s4.60}]
Using the notation of \S~\ref{s4.61} and \S~\ref{s4.62}, replacing $C$ 
by $Y$, we have the following commutative diagram
\begin{equation*}
\xymatrix{
W_{\mathcal{X}}\setminus N_{\mathcal{X}}\ar[rr]\ar[rd]&&
Y^{(n)}\setminus \Delta\ar[dl]^{u}\\
&Pic^{n}Y
}
\end{equation*}
The map $W_{\mathcal{X}}\setminus N_{\mathcal{X}}\to 
S=Pic^{n}Y$ is topologically trivial by \cite{TSSM} Theorem~5.2. The 
map $W_{\mathcal{X}}\setminus N_{\mathcal{X}}\to Y^{(n)}\setminus 
\Delta$ given by $(x,D)\mapsto D$ is a topological unramified covering of 
degree $n$, hence a 
Serre fibration. By Lemma~\ref{s4.60a} we obtain that 
$u:Y^{(n)}\setminus\Delta\to Pic^{n}Y$ is a Serre fibration. The 
canonical exact sequence of homotopy groups and the
equality $\pi_{2}(Pic^{n}Y)=0$ yields \eqref{es4.60}.
\end{proof}

\begin{acknowledgments}
This work was supported 
by a research grant from the University of Palermo.  The author was on leave of absence from the Institute of Mathematics and Informatics of the Bulgarian Academy of Sciences.
\end{acknowledgments}


\end{document}